\documentclass[12pt,reqno]{amsart}
\usepackage{amsmath, amssymb,graphicx,caption, subcaption, mathtools, enumitem}
\usepackage{tikz}
\usepackage{verbatim}
\usepackage{etoolbox}
\patchcmd{\section}{\scshape}{\bfseries}{}{}
\patchcmd{\subsection}{\bfseries}{\itshape}{}{}

\usepackage{color}

\makeatletter
\def\@seccntformat#1{%
  \protect\textup{\protect\@secnumfont
    \ifnum\pdfstrcmp{section}{#1}=0 \bfseries\fi
    \ifnum\pdfstrcmp{subsection}{#1}=0 \itshape \fi
    \csname the#1\endcsname
    \protect\@secnumpunct
  }%
}  
\makeatother

\usepackage[utf8]{inputenc}
\usepackage[english]{babel}

\DeclarePairedDelimiter\norm{\lVert}{\rVert}%
\theoremstyle{plain}

\newtheorem{theorem}{Theorem}[section]
\newtheorem{proposition}[theorem]{Proposition}

\newtheorem{lemma}[theorem]{Lemma}

\theoremstyle{definition}
\newtheorem{definition}[theorem]{Definition}

\newtheoremstyle{note}
{3pt}
{3pt}
{\itshape}
{}
{\itshape}
{:}
{.5em}
{}

\theoremstyle{note}

\newtheorem{example}[theorem]{Example}

\theoremstyle{plain} 
\newcommand{\thistheoremname}{}
\newtheorem*{genericthm}{\thistheoremname}

\DeclareMathOperator{\Id}{Id}

\DeclareMathOperator{\Conv}{Conv}
\DeclareMathOperator{\Ker}{Ker}
\DeclareMathOperator{\Image}{Im}

\DeclareMathOperator{\dist}{dist}

\numberwithin{equation}{section}

\newcommand{\acknowledge}{\subsection*{Acknowledgements}}

\def\ignore#1{}
\def\eatspace#1{\ifx#1 \eatspace\else#1\fi}

\def\x{\mathbf{x}} \def\n{\mathbf{n}}   \def\y{\mathbf{y}} \def\z{\mathbf{z}} 
\def\X{X} \def\Q{Q} \def\I{I}
\def\nv{\kappa }

\def\Ue{U_\varepsilon }

\def\ff{\mathbf{f}} 
\def\ffe{\ff_\varepsilon }

\def\bitems#1{{\parskip=-10pt\parindent=0pt\def\\{\hfill\par\hangindent60pt\ostrut11\qquad $\bullet$ \hskip5pt\eatspace}\\ #1\par}\smallskip}

\def\is#1{{\em #1}} 

\def\rbox#1{\hbox{\rm #1}}
\def\roh{\rbox}

\def\rox#1{\quad \rbox{#1}\quad }
\def\roy#1{\ \ \rbox{#1}\ \ }
\def\roq#1{\qquad \rbox{#1}\qquad }

\def\ldotsx{\ \ldots\ } \def\cdotsx{\ \cdots\ }

\def\R{\mathbb{R}}  \def\Rd{\R^d}
\def\Sm{\mathbb{S}} \def\Sdo{\Sm^{d-1}}  \def\So{\Sm^1}

\def\barr#1{\overline{#1}{}}

        \let\bpa\bparens

     \let\bbk\bbrackets

  \let\bco\braces

\def\cupz{\,\cup\,} \def\capz{\,\cap\,}

\def\setb#1{\{#1\}} 
\def\szero{\setb0} 

\def\hexnumber#1{\ifcase#1 0\or1\or2\or3\or4\or5\or6\or7\or8\or9\or
 A\or B\or C\or D\or E\or F\fi}

\edef\msbhx{\hexnumber\symAMSb}   
\mathchardef\emptyset="0\msbhx3F

\def\operator#1{\expandafter\def\csname#1\endcsname{\mathop{\rm #1}\nolimits}}

\operator{sign} \operator{mod}

\def\O{{\rm O}}

\def\Labl#1{{\it#1\/}:\enspace}
\def\Label#1{\medbreak\noindent\Labl{#1}}
\def\Labelpar #1.#2\par{\Label{#1}#2\par\medbreak}

\def\Proof{\Label{Proof}}

\def\qed{\hfill{\it Q.E.D.}\par}

\def\ophantom#1#2{\setbox0=\hbox{$#1#2$}\setbox2 = \null
                 \ht2 = \ht0 \dp2 = \dp0 \box2}
\def\odphantom#1{\ophantom\displaystyle{#1}}
\def\otphantom#1{\ophantom\textstyle{#1}}

\def\set#1#2{\mathchoice{\left \{\,#1\odphantom{#2}\;\right |
\left .\;#2\odphantom{#1}\,\right \}}{\{\,#1\otphantom{#2}\,|\,#2\otphantom{#1}\,\}}
{\{#1\otphantom{#2}\,|\,#2\otphantom{#1} \}}
{\{#1\otphantom{#2}\,|\,#2\otphantom{#1} \}}}

\def\ip#1#2{\langle \,{#1}\,,{#2}\,\rangle}

\def\Limu#1{\lim_{#1}\>} \def\Lim#1#2{\Limu{#1\to #2}}

  \def\forally{\roy{for all}}

\def\rgs#1#2#3{#1=#2,\ldots,#3} 
\def\rg#1#2{#1=1,\ldots,#2} \def\rgo#1#2{#1=0,\ldots,#2} 

\def\subs #1#2{#1_1,\ldots,#1_{#2}}

\def\psubos #1#2{(#1_0,\ldots,#1_{#2})}
\def\bsubs #1#2{\{#1_1,\ldots,#1_{#2}\}}

\def\Summ#1#2#3{\sum _{#1\:=\:#2}^{#3}\>}
\def\Sum#1#2{\Summ{#1}1{#2}}  

\def\comp{\raise 1pt \hbox{$\,\scriptstyle\circ\,$}}

\def\longrightarrowx{\ \longrightarrow\ } 
\def\longrightarrowz{\,\longrightarrow\,}

\def\simarrow{\setbox0 = \hbox{$\longrightarrow$} 
\setbox2 = \hbox to \wd0{\hfil$\widetilde {\hbox to .8 \wd0{\hfill}}\,$\hfil}
\wd2 = \wd0 \wd0=0pt \; \box0 \lower 1pt \box2\>}

\def\Eq#1$$#2$${\StEq#1  \EnEq{#2}}
\def\StEq#1 #2\EnEq#3{\begin{equation}\label{#1eq*} #3\end{equation}}

\def\eqrefz#1{\ref{#1eq*}}
\def\eq#1{{\rm (\eqrefz{#1})}}
\def\eqs#1#2{\eq{#1},~\eq{#2}}
\def\eqss#1#2#3{\eq{#1},~\eq{#2},~\eq{#3}}

\newdimen\eqjot \eqjot = 1\jot
\def\openupeq{\openup \the\eqjot}

\def\qaeq#1#2{{\def\\{&}\vcenter{\openupeq\halign{$\displaystyle
   ##\hfil$&&\hskip#1pt$\displaystyle##\hfil$\cr #2\cr}}}}
 
\def\req{\qaeq{30}}  
\def\qeq{\qaeq{20}} \def\weq{\qaeq{15}}
\def\xeq{\qaeq{10}}  

\def\addtab#1={#1\;&=}
\def\addtabe#1=#2={#1=#2\;&=}

\def\ezeq#1#2#3{{\def\\{\cr#1}\vcenter{\openupeq \halign{$\displaystyle 
   \hfil##$&$\displaystyle##\hfil$&&\hskip#2pt$\displaystyle##\hfil$\cr#1#3\cr}}}}
\def\eaeq{\ezeq\addtab}

\def\eeq{\eaeq{20}}

\def\saeq#1#2{{\def\\{\cr}\vcenter{\openup1\jot \halign{$\displaystyle
   ##\hfil$&&\hskip#1pt$\displaystyle##\hfil$\cr #2\cr}}}}

\def\macases#1#2{\left\{\enspace\saeq{#1}{#2}\right.}
\def\mcases{\macases{20}}
\def\xcases{\macases{15}}

\def\ostrut#1#2{\hbox{\vrule height #1pt depth #2pt width 0pt}}
\def\sstrut{\ostrut0}

\begin{document}

\title[Maps Converging to Convex Hulls]{Continuous Maps from Spheres Converging to Boundaries of Convex Hulls}
\author{J.~Malkoun}
\address[J.~Malkoun]{Department of Mathematics and Statistics\\
Faculty of Natural and Applied Sciences\\
Notre Dame University-Louaize, Zouk Mosbeh\\
P.O.~Box: 72, Zouk Mikael, Lebanon}
\email[J.~Malkoun]{joseph.malkoun@ndu.edu.lb}
\author{P.J.~Olver}
\address[P.J.~Olver]{School of Mathematics\\
University of Minnesota\\
Minneapolis, MN   55455, USA}
\email[P.J.~Olver]{olver@umn.edu}
\date{\today}

\begin{abstract} Given $n$ distinct points $\mathbf{x}_1, \ldots, \mathbf{x}_n$ in $\mathbb{R}^d$, let $K$ denote their convex hull, which we assume to be $d$-dimensional, and $B = \partial K $ its $(d-1)$-dimensional boundary.  We construct an explicit one-parameter family of continuous maps $\mathbf{f}_{\varepsilon} \colon \mathbb{S}^{d-1} \to K$ which, for $\varepsilon > 0$, are defined on the $(d-1)$-dimensional sphere and have the property that the images $\mathbf{f}_{\varepsilon}(\Sdo)$ are codimension $1$ submanifolds contained in the interior of $K$.  Moreover, as the parameter $\varepsilon$ goes to $0^+$, the images $\mathbf{f}_{\varepsilon}(\Sdo)$ converge, as sets, to the boundary $B$ of the convex hull. We prove this theorem using techniques from convex geometry of (spherical) polytopes and set-valued homology. We further establish an interesting relationship with the Gauss map of the polytope $B$, appropriately defined. Several computer plots illustrating our results will be presented.
\end{abstract}

\maketitle

\section{Introduction}

Given a configuration $\X = (\x_1,\ldots,\x_n)$ of $n$ distinct points in $\R^d$, computing their convex hull $K = \Conv(\X)$ is a 
famous problem in Computational Geometry. Many algorithms have been developed for  this task, including the Gift Wrap or Jarvis March algorithm, the Graham Scan algorithm, QuickHull, Divide and Conquer, Monotone Chain or Andrew's algorithm, Chan's algorithm, the Incremental Convex Hull algorithm, the Ultimate Planar Convex Hull algorithm, and others.  See, for instance, \cite{CG} and the references within.  

In this paper, we develop an alternative, direct approach to this problem that does not rely on any underlying computer algorithm. Instead, assuming $\dim K = d$, meaning that its interior $K^\circ$ is a nonempty open subset of $\Rd$, we construct a one-parameter family of approximations to its $(d-1)$-dimensional boundary $B = \partial K$, as the images of continuous maps $\ffe\colon \Sdo \to \R^d$ for $ \varepsilon >0 $, that are defined explicitly, and fairly simply, in terms of the points $\x_1,\ldots,\x_n$.

Initial computer generated plots suggested that the images $\ffe( \Sdo)$ of our family of maps provide excellent approximations to the boundary $B$ for all configurations that we have tried; see Figures \ref{fig1} and \ref{fig2} for some representative examples. Our main result, Theorem \ref{main_theorem}, states that the images $\ffe(\Sdo)$  converge, as sets, to the boundary $B$ as the parameter $\varepsilon \to 0^+$.  We will  also explain in detail the mechanism of convergence. We then establish a relationship with the Gauss map of a smooth surface, \cite{GAS}, thereby defining the inverse Gauss map of the boundary of the convex hull as a set-valued map.  Indeed, our proof of the main theorem relies on techniques from the theory of set-valued homology.

On the other hand, the convergence of the approximating sets $\ffe (\Sdo)$ to the boundary $B$ is highly non-uniform.  Indeed, as we will see, the images $\ffe(\n)$ of almost every point $\n\in \Sdo$ converge to one of the vertices of $B$.  Thus, if one discretely samples $\Sdo$ by a large but finite number of points $\subs \y N$, most of their image points $\ffe(\y_k) \in \ffe(\Sdo)$ will accumulate around the vertices of $B$, and the remainder of $B$ will be increasingly sparsely approximated as $\varepsilon  \to 0^+$.  This non-uniform sampling property can be observed in the three-dimensional illustrative plots in Figure \ref{fig2}.  



See \cite{Firey} for an alternative, less practical approach to approximating convex polytopes and convex sets by algebraic sets. A potential future project based on these constructions will be to develop fast practical algorithms for approximating the convex hull of a point configuration.  A potentially interesting extension of our techniques will be to the approximation of Wulff shapes of crystals, \cite{NS}.

\medskip

\section{A Family of Maps defined by a Point Configuration} \label{def}

Let us begin by introducing the basic set up and our notation, before defining the family of maps that will be our primary object of study. 

Let $C_n(\R^d)$ denote the configuration space of $n$ distinct points in $\R^d$. Let $\X = (\x_1,\ldots,\x_n) \in C_n(\R^d)$, 
so that each $\x_i \in \R^d$ and $\x_i \ne \x_j$ whenever $i\ne j$. Assuming $n \geq d+1$, let $C_n^*(\R^d) \subset C_n(\R^d)$ denote the dense open subset of \emph{nondegenerate configurations}, meaning those whose points do not all lie on a proper affine subspace of $\Rd$.  From here on we fix the nondegenerate point configuration $\X \in C_n^*(\R^d)$, and suppress all dependencies thereon. 

Let $K = \Conv(\X) \subset \Rd$ denote the convex hull of the points in $X$, which, by nondegeneracy, is a bounded convex polytope of dimension $d$ whose interior is a nonempty open subset $K^\circ \subset \Rd$, \cite{Grunbaum,Ziegler}.  Let $B = \partial K = \partial \Conv(\X)$ be its boundary, which is a piecewise linear closed hypersurface in $\Rd$, forming a $(d-1)$-dimensional polytope.

Given any pair of indices $1 \leq i,j \leq n$ with $i \neq j$, we define real-valued functions
$c_{ij}\colon  \R^+ \times \Sdo \to \R^+ = \{0 < t \in \R\}$ 
by
\Eq{cij}
$$\weq{c_{ij}(\varepsilon, \n)  =\varepsilon + \max\{ 0, \> -\,\ip{\n}{\n_{ij}} \},\\\   \varepsilon > 0, \\ \n \in \Sdo,}
$$
where $\ip\cdot\cdot$ denotes the Euclidean inner product in $\mathbb{R}^d$, and where
\Eq{nij}
$$\req{\n_{ij}  = \frac{\x_j - \x_i}{\norm{\x_j - \x_i}} \in \Sdo,\\ i\ne j,}$$
is the unit vector pointing from $\x_i$ to $\x_j$, with $\norm {\,\cdot\,}$ denoting the Euclidean norm.  Note that $\n_{ij}  = -\,\n_{ji}$.  The $c_{ij}$ in \eq{cij} are continuous maps; moreover, $c_{ij}(\varepsilon,\n) > 0$ since we are assuming (for now) that $\varepsilon >0$. We further define, for any $1 \leq i \leq n$, the map $c_i\colon \R^+ \times \Sdo \to \R^+$ by the $(n-1)$-fold product
\Eq{ci}
$$ c_i(\varepsilon,\n) = \prod_{\substack{1 \leq j \leq n \\ j \neq i}} c_{ij}(\varepsilon,\n). $$
Finally, let us set
\Eq{li}
$$\req{\lambda _i(\varepsilon,\n) = \frac{c_i(\varepsilon,\n)}{\Delta (\varepsilon,\n)},\\ \rg in,}$$
where
\Eq{Delta}
$$\weq{
\Delta (\varepsilon,\n) = \sum_{j=1}^n c_j(\varepsilon,\n) > 0 \rox{for all} \varepsilon > 0,\\ \n \in \Sdo.}$$
Thus,
\Eq{ls}
$$\qeq{0 < \lambda _i(\varepsilon,\n) < 1, \rox{and} \ \ 
\Sum in\lambda _i(\varepsilon,\n) = 1.}$$

Given a point configuration $\X \in C_n(\R^d)$, we can now define the main object of interest in this paper: the one-parameter family of maps
$ \ffe\colon  \Sdo \to \R^d $ defined by
\Eq{fe}
$$\weq{\ffe (\n) = \Sum in\lambda _i(\varepsilon,\n) \x_i,\\ \ \ \varepsilon > 0,\\ \n \in \Sdo.}$$
From \eq{ls}, \eq{fe}, one immediately deduces that
$$\ffe(\n) \in K^\circ, \rox{for any} \weq{\varepsilon > 0,\\\n \in \Sdo.} $$

\begin{figure} 
\begin{tabular}{ccc}
 \hskip-.5in \includegraphics[width=60mm,trim=0mm 85mm 0mm 70mm]{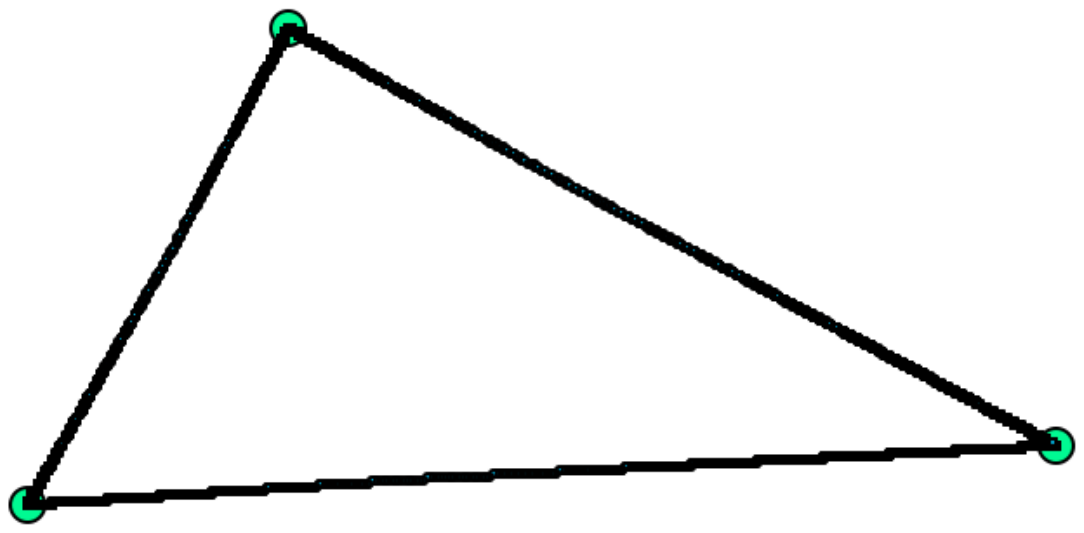} &
 \hskip-.85in   \includegraphics[width=60mm,trim=0mm 85mm 0mm 70mm]{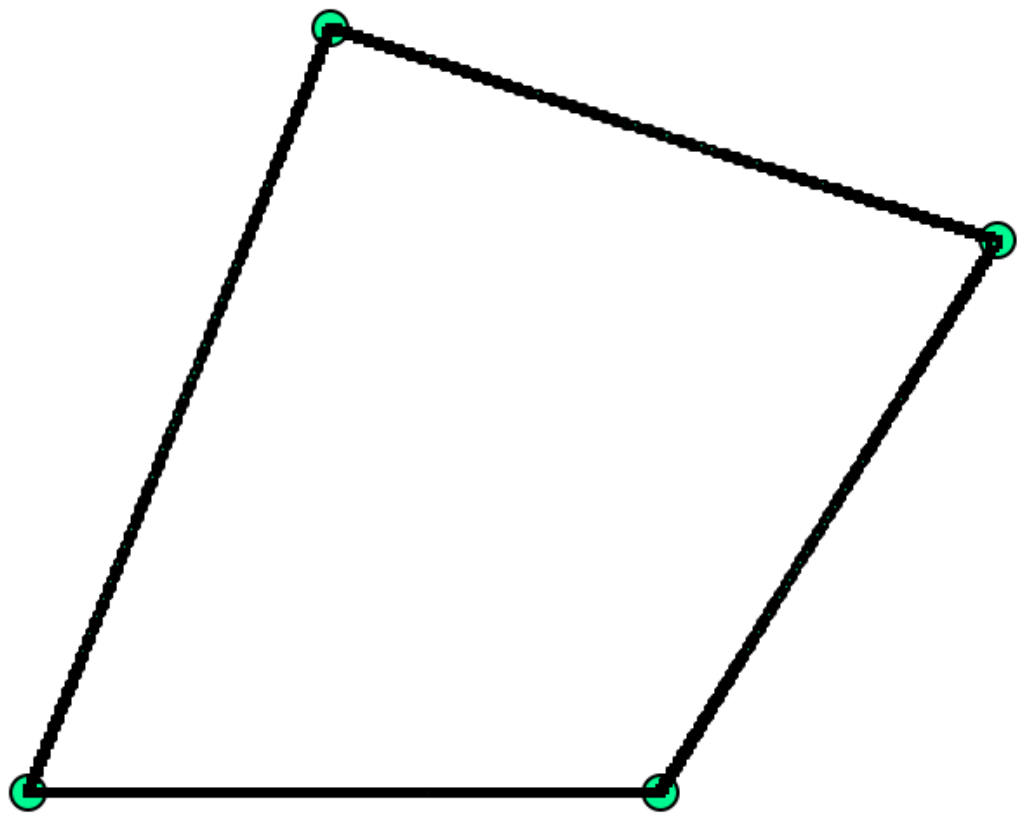} 
 &\hskip-.85in  \includegraphics[width=60mm,trim=0mm 85mm 0mm 70mm]{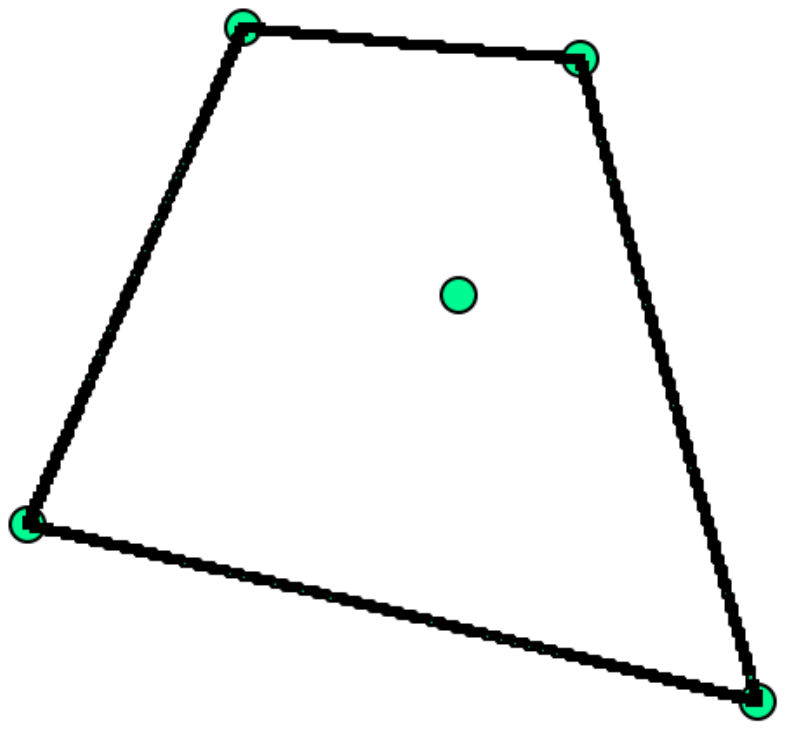}\\
\hskip-.5in(a) $n=3$ &\hskip-.85in (b) $n=4$ &\hskip-.85in (c) $n=5$ \\[6pt]
\hskip-.5in\includegraphics[width=60mm,trim=0mm 85mm 0mm 70mm]{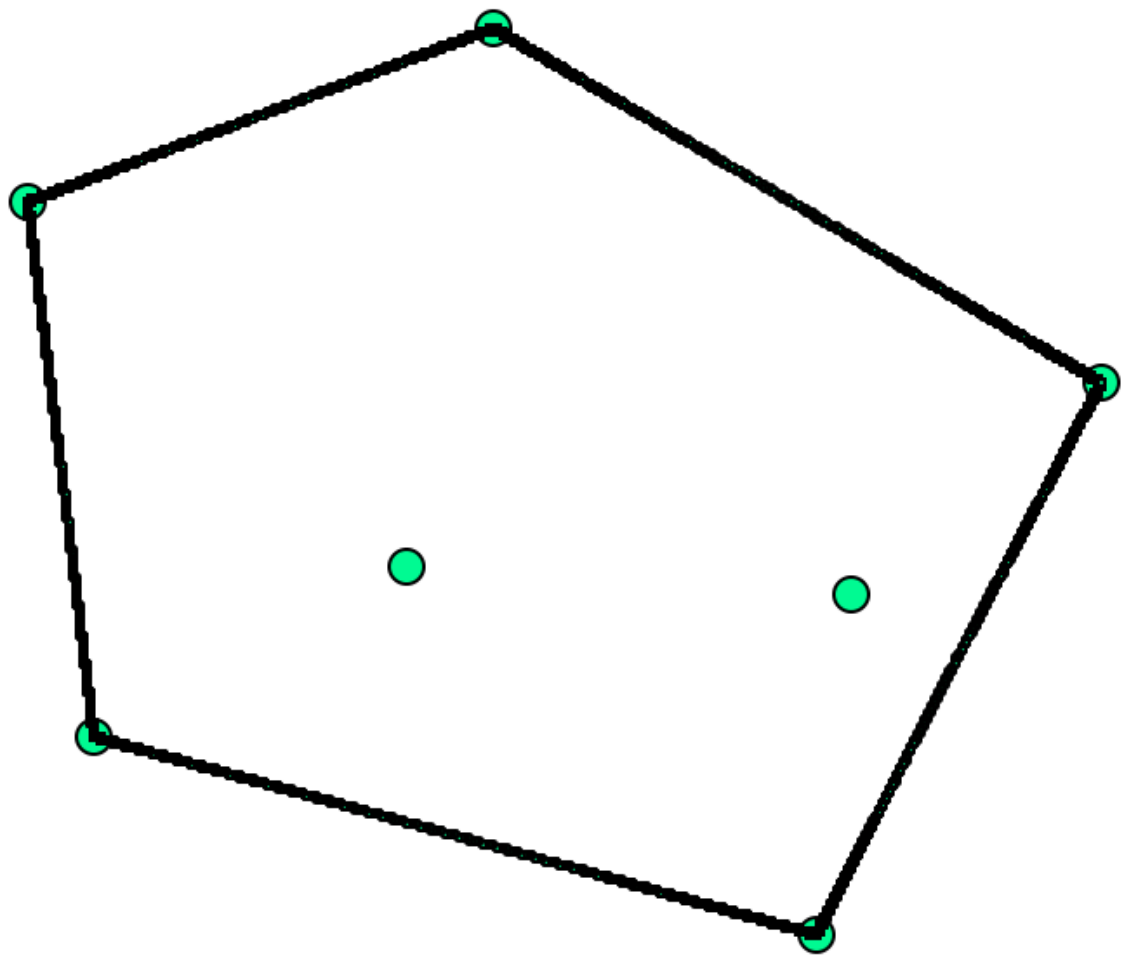} &  
\hskip-.85in \includegraphics[width=60mm,trim=0mm 85mm 0mm 70mm]{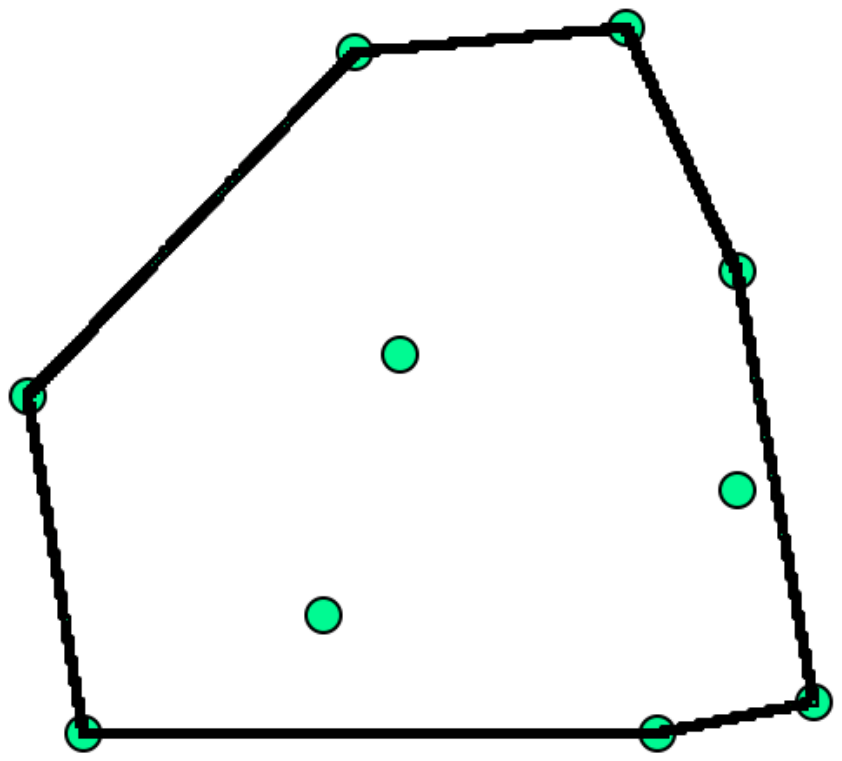} & 
\hskip-.85in \includegraphics[width=60mm,trim=0mm 85mm 0mm 70mm]{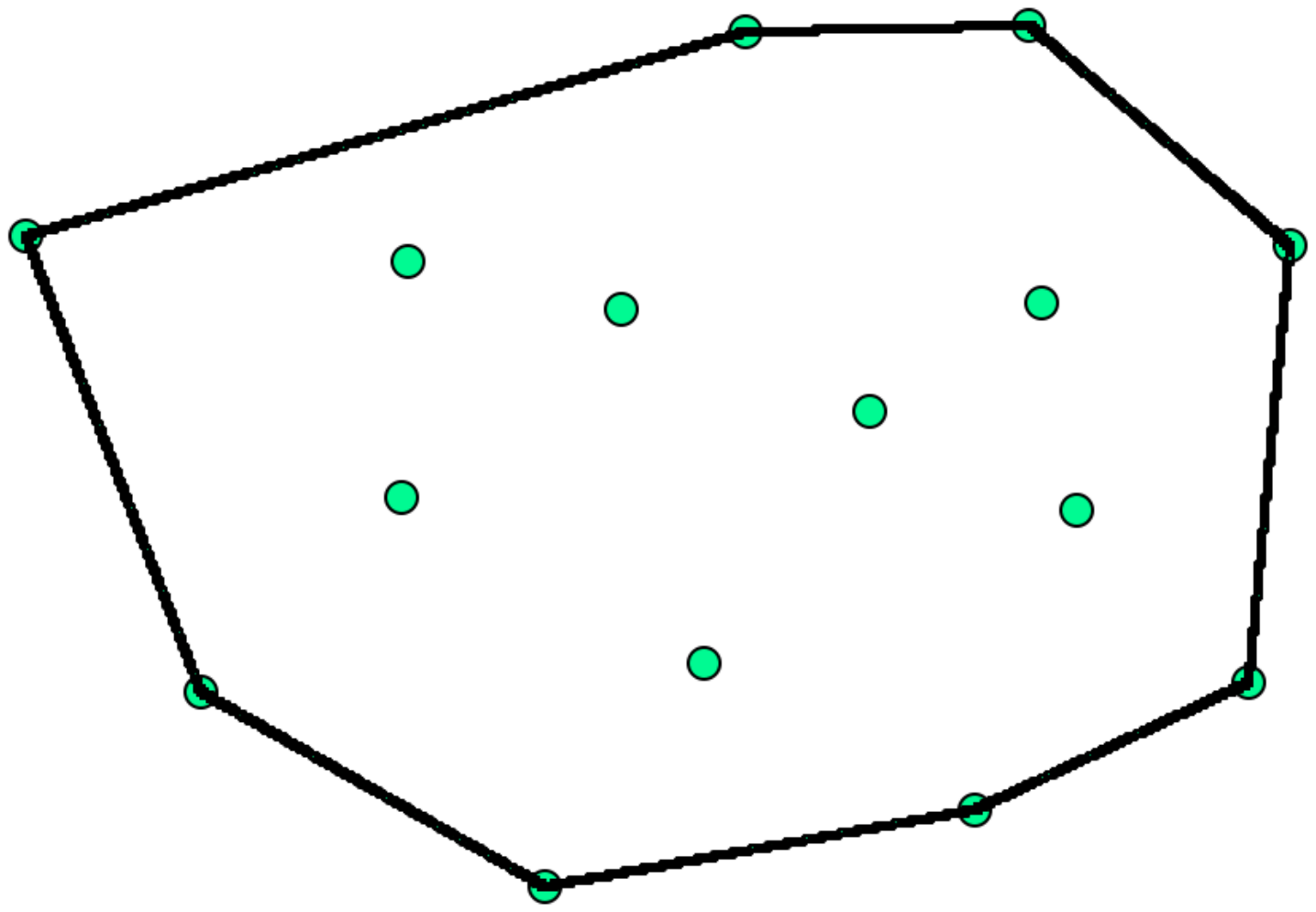}\\
\hskip-.5in(d) $n=7$ & \hskip-.85in(e) $n=10$ & \hskip-.85in(f) $n=15$ \\[6pt]
\end{tabular}
\caption{Plots of planar point configurations $\X$ and the image $\ffe(\So)$, with $\varepsilon = 0.01$.} 
\label{fig1}
\end{figure}

\begin{figure} 
\begin{tabular}{c}
\includegraphics[width=100mm,trim=0mm 75mm 0mm 70mm]{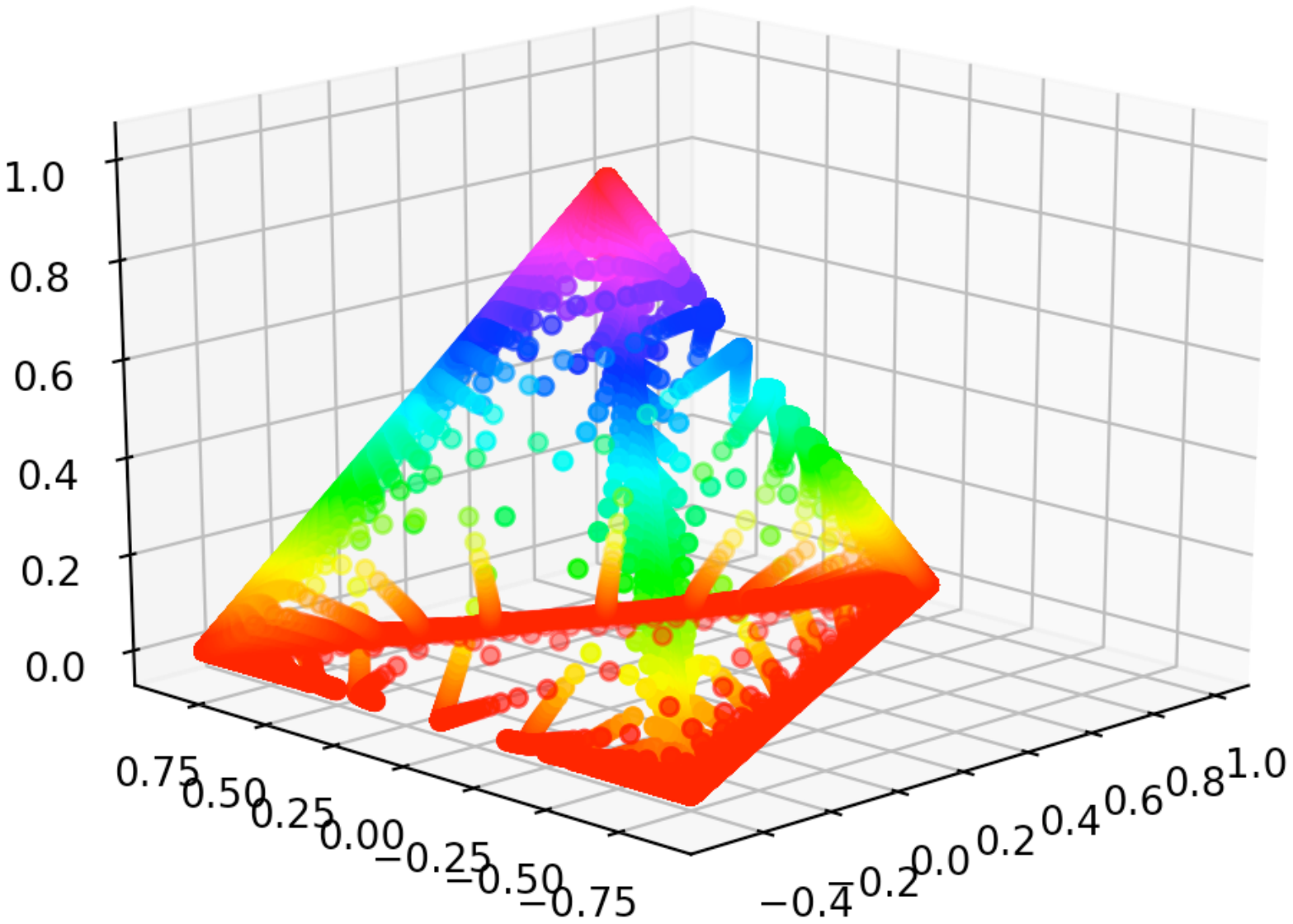} \\
(a) $\X$ consists of the vertices of a regular tetrahedron, so $n=4$. \\ 
\includegraphics[width=100mm,trim=0mm 70mm 0mm 70mm]{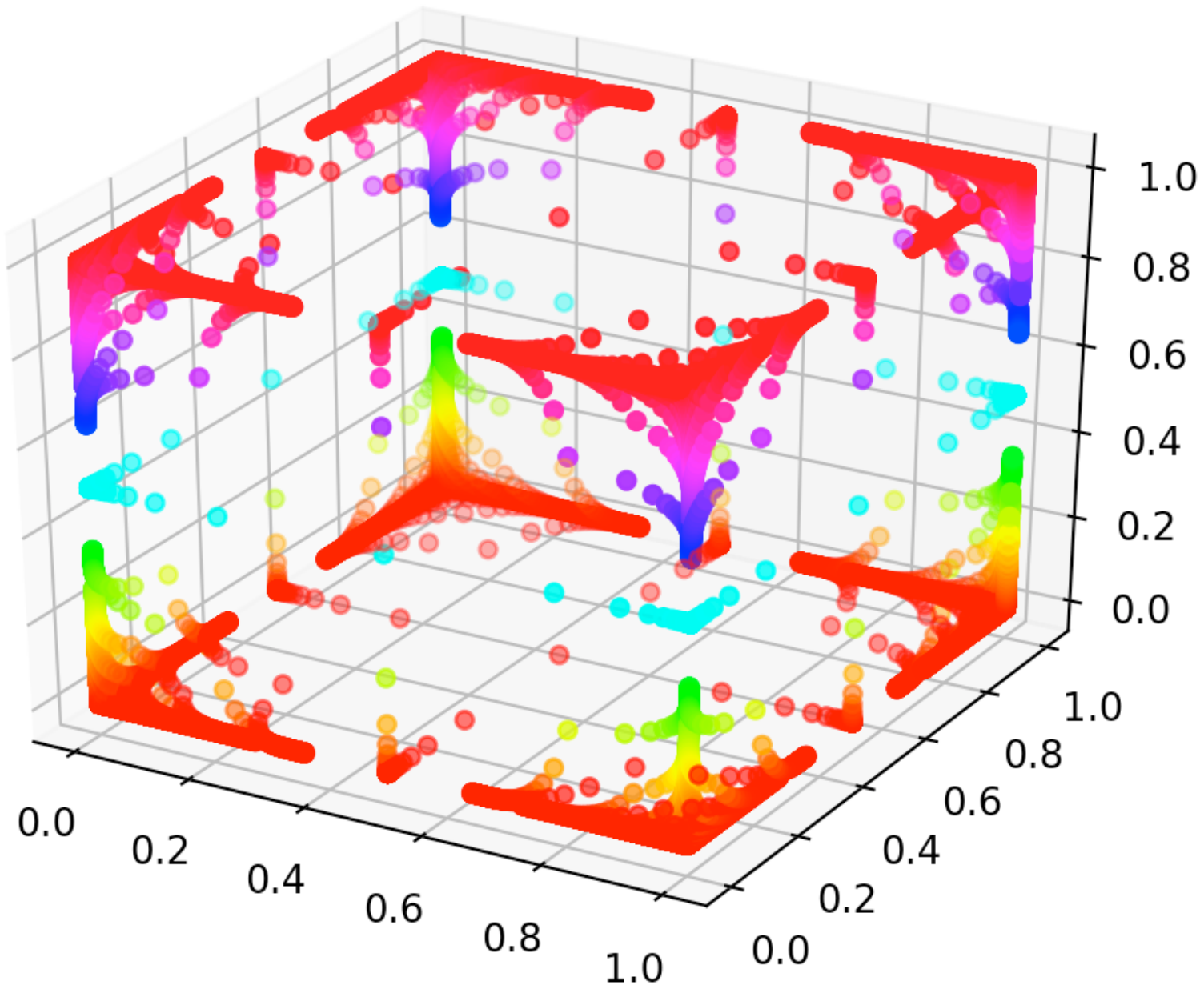} \\
(b) $\X$ consists of the vertices of a cube, so $n=8$.\\[6pt]
\end{tabular}
\caption{Plots of point configurations $\X$ in dimension $d=3$ and the images of sample points on $\Sm^2$ under the map $\ffe$ with $\varepsilon = 0.01$.}
\label{fig2}
\end{figure}

\begin{figure} 
\centering
\includegraphics[width=60mm,trim=0mm 0mm 0mm 0mm]{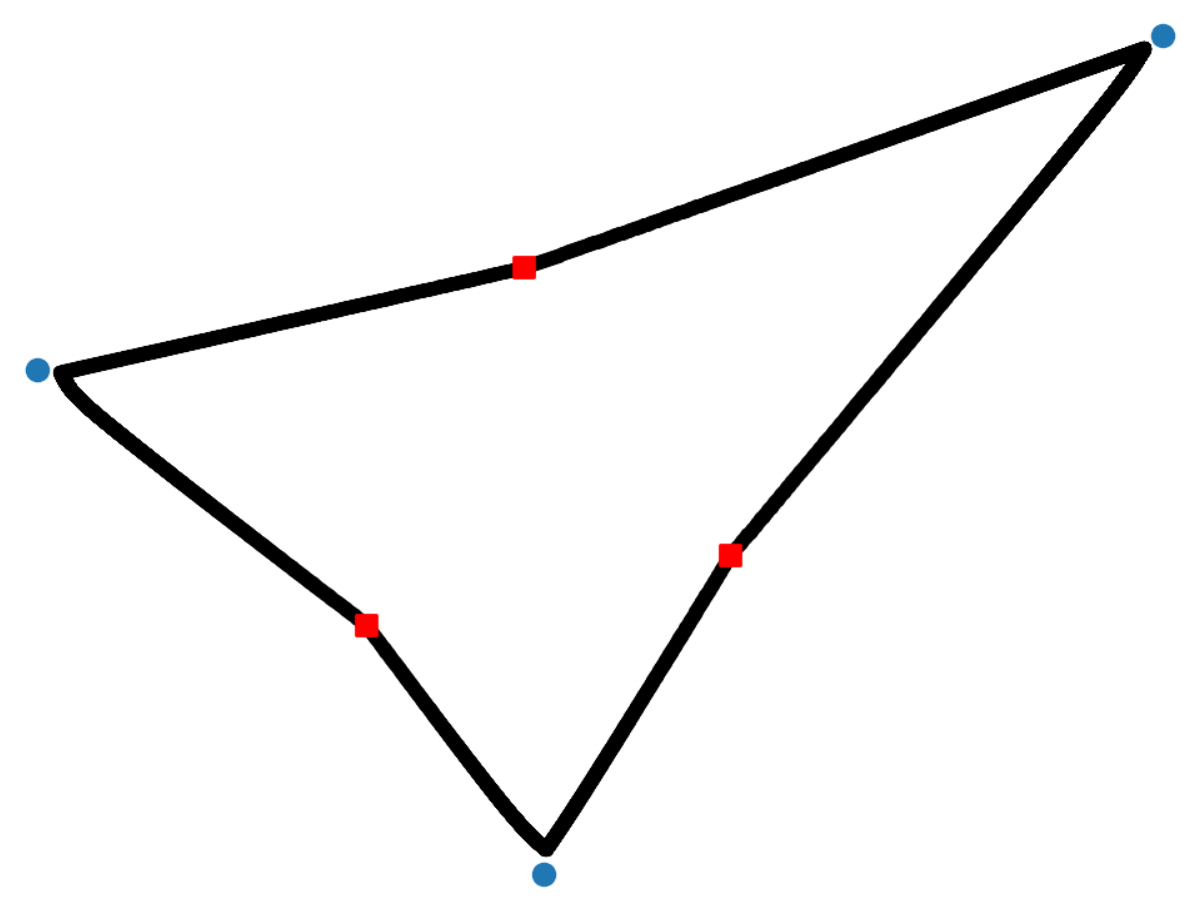} \\
\caption{$X$ consists of the vertices of a triangle, with $\varepsilon = 0.1$. One can see that $\ffe(\Sm^1)$ is indented at the points $\ffe (\n)$, shown as small red rectangles, which correspond to 
$\n \in S_{ij}$, for $1 \leq i,j \leq 3$ and $i \neq j$.}
\label{fig3}
\end{figure}

Inspection of Figures \ref{fig1} and \ref{fig2}, and others that can be easily generated by computer, indicates that, for a given $\X \in C^*_n(\R^d)$ and small $\varepsilon > 0$, the image of $\Sdo$ under $\ffe$ may be 
used as a good approximation of the boundary $B = \partial K \subset \Rd$ of the convex hull of $\X$. More precisely, the Main Theorem to be proved in this paper is as follows:

\begin{theorem} \label{main_theorem} Given $\X \in C_n^*(\R^d)$, let $K = \Conv(\X)$ be their convex hull, which has dimension $d$.  Let $\ffe$ be defined by \eq{fe}.  Then, for $\varepsilon > 0$, the images of the unit sphere under $\ffe$ lie in the interior of the convex hull of $\X$, so  $\ffe(\Sdo) \subset K^\circ$, and, moreover, converge to its boundary as sets in $\R^d$\/{\rm:} 
\Eq{main}
$$\Lim\varepsilon{0^+} \ffe(\Sdo) = \partial K.$$
\end{theorem}

The set theoretic convergence in \eq{main} is \emph{uniform} in the sense that the images $\ffe(\Sdo)$ lie in an $\O(\varepsilon )$ neighborhood of the boundary $\partial K$, even though their pointwise convergence is highly nonuniform.  See below for precise details on what this means.

\medskip

\emph{Remark\/}:  On the other hand, if the point configuration is degenerate, meaning $\X \in C_n(\R^d) \setminus C_n^*(\R^d)$ and so its convex hull $K = \Conv(\X)$ has dimension $< d$, then one can show that $\Lim\varepsilon{0^+} \ffe(\Sdo) = K$.  Indeed, observe that the maps $\ffe$ depend continuously on the point configuration.  If one slightly perturbs $\X$ to a nondegenerate configuration $\X_\delta \in C_n^*(\R^d)$, then their perturbed convex hull $K_\delta $ is of dimension $d$ and, by the Theorem, $\ffe(\Sdo) \to \partial K_\delta $.  But as $\delta \to 0$, their boundaries converge to the entire convex hull: $\partial K_\delta \to  K$, which enables one to establish the result.  Since this case is of less importance for our purposes, the details are left to the reader.

\medskip

%

In Section \ref{nc}, we present notions from convex geometry that are relevant to this work, including normal cones and normal spherical polytopes.  The latter enable us to associate with a convex polytope $B$ a spherical complex $S^*_B$, \cite{MS}, meaning a tiling of $\Sdo$ by spherical polytopes, with the property that it has the same combinatorial type as the dual polytope $B^*$. Then, in Section \ref{Gauss}, we connect our 
constructions with the differential geometric  concept of the Gauss map of a convex hypersurface, generalized to  the boundary of the convex polytope.  We explain how our maps converge to the inverse Gauss map  of the boundary of the convex hull of the point configuration, which is viewed as a set-valued function.  Finally in Section 
\ref{main}, we prove our main result using a combination of convex geometry and set-valued homology theory, the latter described in Appendix \ref{appB}.


\section{Convex Geometry and (Spherical) Polytopes} \label{nc}


Let us recall some basic terminology and facts about convex sets and cones, and both flat and spherical polytopes, many of which can be found in \cite{BV, Rockafellar}. The closed cones appearing in this paper are convex, \emph{pointed}, meaning they do not contain any positive dimensional linear  subspace of $\Rd$, and \emph{polyhedral}, meaning they can be characterized as the intersection of finitely many, and at least two, closed half spaces, \cite{CP,Ziegler}. On the other hand, for us an \emph{open cone} $N \subset \Rd$ is a cone such that $N \setminus \{0\}$ is an open subset of $\Rd$ and such that its closure $\barr N$ is of the above type.

Let us fix a nondegenerate point configuration $\X \in C_n^*(\R^d)$ consisting of $n$ distinct points $\subs \x n \in \Rd$.     Let $K = \Conv(\X) \subset \Rd$ denote the convex hull of the points in $X$, which is a bounded convex $d$-dimensional polytope, \cite{Grunbaum,Ziegler}. Let $B = \partial K = \partial \Conv(\X)$ be its boundary, which is itself a polytope of dimension $d-1$ --- specifically a 
piecewise linear closed hypersurface in $\Rd$.   Assume, by relabelling if necessary, that $\x_1,\ldots,\x_\nv$ are the vertices of $K$, while 
$\x_{\nv+1}, \ldots, \x_n$ are the remaining points, which may either lie in the interior $K^\circ$ or at a non-vertex point of the boundary $B$.   
 The \is{faces} of $B$ range in dimension from $0$, the \is{vertices}, to $1$, the \is{edges}, up to $d-1$, the \is{facets}. Two vertices are \is{adjacent} if they are the endpoints of a common edge. Note that each face $F \subset K$ is itself a convex polytope. If $0 < m \leq d-1$, we denote the \emph{interior} of an $m$-dimensional face $F$ by $F^\circ = F \setminus \partial F$, which is a flat $m$-dimensional submanifold of $\Rd$.  (Keep in mind that this is not the same as its interior as a subset of $\Rd$, which is empty.)  

Define the \is{normal cone} at the point $\x_i$ by
\Eq{Ni}
$$N_i = \set{\y \in \Rd}{\ip \y{\n_{ij}} \leq 0 \forally j\ne i} = \bigcap_{j \ne i} \>H_{ij}, $$
where the  unit vectors $\n_{ij} \in \Sdo$ are given in \eq{nij}, and
\Eq{Hi}
$$H_{ij} = \set{\y \in \Rd}{\ip \y{\n_{ij}} \leq  0} $$
is the closed half space opposite to $\n_{ij}$.
 Further let
\Eq{Nio}
$$N^\circ_i = \set{\y \in \Rd}{\ip \y{\n_{ij}} < 0 \forally j\ne i}.$$
denote the interior of the normal cone $N_i$.
It is easy to see that $N^\circ_i \ne \emptyset$ if and only if $\x_i$ is a vertex.   Also, $N^\circ_i \cap N^\circ_j = \emptyset $ whenever $i \ne j$.  Indeed, if $\y \in N^\circ_i$ then $\ip \y{\n_{ij}} < 0$. But then $\ip \y{\n_{ji}} = \ip \y{-\n_{ij}} > 0$ and hence $\y \not\in N^\circ_j$.  
Furthermore, the union of the vertex normal cones is the entire space:
\Eq{Nunion}
$$\bigcup_{i=1}^\nv \>N_i = \Rd,$$
i.e., every vector is in one of the normal cones. This is a direct consequence of the Supporting Hyperplane Theorem; see for instance \cite[pp. 50--51]{BV}. 


%

A \is{spherical polytope} is characterized as the intersection of finitely many closed hemispheres that does not contain any antipodal points, cf.~\cite[\S 2.2]{CP}. It can alternatively be characterized as the intersection $C \capz \Sdo$ of the unit sphere with a pointed polyhedral cone $C \subset \Rd$.
Let us consequently define the \is{normal spherical polytope}
\Eq{Si}
$$S_i = N_i \capz \Sdo  = \set{\n \in \Sdo }{\ip \n{\n_{ij}} \leq  0 \forally j\ne i},$$
associated with the point $\x_i$. Its interior
\Eq{Sio}
$$S^\circ_i = N^\circ_i \capz \Sdo  = \set{\n \in \Sdo }{\ip \n{\n_{ij}} <  0 \forally j\ne i}$$
 is  nonempty if and only if $\x_i$ is a vertex, in which case it is an open submanifold of the unit sphere.  Note that, by \eq{Nunion} and the preceding remarks,
\Eq{SiS}
$$\qeq{\bigcup_{i=1}^\nv \>S_i = \Sdo ,\\S^\circ_i \cap S^\circ_j = \emptyset ,\\i \ne j.}$$

The \is{normal cone} and \is{normal spherical polytope} associated with a general point $\x \in K$ in the convex hull are similarly defined:
\Eq{NS}
$$\eeq{N_\x = \set{\y \in \Rd}{\ip \y{\z-\x} \leq 0 \forally \z \in K},\\S_\x = N_\x \capz \Sdo  = \set{\n \in \Sdo }{\ip \n{\z-\x} \leq 0 \forally \z \in K}.}$$
As above, $N_\x = \szero$ if $\x \in K^\circ$, while $N_{\x_i} = N_i$ when $\x_i$ is a vertex. More generally, if $F \subset B = \partial K$ is an $m$-dimensional face, then the normal cone $N_\x$ is independent of the point $\x \in F^\circ$ in its interior, and we thus define $N_F = N_\x$ for any such $\x \in F^\circ$.  If the face $F$ has dimension  $m$, then $N_F$ is a $(d-m)$-dimensional cone. Define its \emph{interior} to be $N^\circ_F = N_F \setminus \partial N_F$, which is a $(d-m)$-dimensional submanifold of $\Rd$. {\it Warning\/}: unless $m=0$, so $F$ is a vertex, $N^\circ_F$ is not the same as the interior of $N_F$ considered as a subset of $\Rd$, which is empty. In particular, if $F$ is a facet, i.e., a $(d-1)$-dimensional face, then $N_F$ is a one-dimensional cone, i.e., a ray in the direction of its unit outward normal $\n_F$, with $N^\circ_F = \set {c\, \n_F}{c> 0}$. Observe that if $H \subset \partial F$ is a subface, then $N_F \subset \partial N_H$.  Further, convexity of $K$ implies that $N_F^\circ \capz N_G^\circ = \emptyset$ whenever $F \ne G$ are distinct faces of $B$; in particular, $\n_F \ne \n_G$ whenever $F \ne G$ are distinct facets.

The collection of the interiors of all the normal cones to the faces of $B$ form the \is{complete normal fan} associated with the polytopes $B$ and $K$, and their disjoint union fills out the entire space, except for the origin (which can be identified with $N_K$):
\Eq{fan}
$$\Rd = \{0\}  \>\vee \  \bigvee_{F \subset B} N_F^\circ.$$


We further define the \emph{normal spherical polytope} associated with the $m$-dimensional face $F$ as $S_F = N_F \capz \Sdo$. When $m < d-1$, its interior $S^\circ_F = N^\circ_F \capz \Sdo$ is a $(d-m-1)$-dimensional submanifold of $\Sdo$, while for $m=d-1$, the normal spherical polytope $S_F$ is a single point, namely the facet's unit outward normal $\n_F$.  As an immediate consequence of the complete normal fan decomposition \eq{fan}, we can write the sphere as a disjoint union
\Eq V
$$\Sdo = \bigvee_{\dim F < d-1} S^\circ_F \ \vee \bigvee_{\dim F = d-1} S_F,$$
where the second term runs over the facets and the first over all other faces of $B = \partial K$.  As above, if $G \subset \partial F$ is a subface, then $S_F \subset \partial S_G$.

The collection of all normal spherical polytopes $S_F$, where $F$ runs over all faces of $B$, forms a spherical complex, \cite{MS}, denoted $S^*_B$, that tiles the sphere by spherical polytopes as shown in \eq V.  We note that $S^*_B$ has the same combinatorial type as the dual polytope  $B^*$, \cite{Grunbaum}, and hence we regard the \emph{normal spherical complex} $S^*_B$ as the \emph{spherical dual} to $B$.

%
%
\medskip

\emph{Remark\/}:  Maehara and Martini, \cite{MM}, propose a similar construction, that they call the ``outer normal transform'' of a convex polytope $B \subset \Rd$  of dimension $d-1$.  They associate each facet $F \subset B$ with its outward normal $\n_F \in \Sdo \subset \Rd$.  The outer normal transform of $B$ is defined to be the convex hull of the facet normals in $\Rd$.  They observe that, unlike our spherical dual, their transform is not necessarily combinatorially equivalent to the dual polytope $B^*$.  

On the other hand, if we \emph{flatten} all the normal spherical polytopes of the spherical dual $S^*_B$, meaning we replace each $S_F \subset \Rd$ by the convex hull of its vertices, the result will be a polytope $\widehat B \subset \Rd$ contained within the unit ball, all of whose vertices lie on the unit sphere.  Although the resulting polytope $\widehat B$ also has the same combinatorial type as  $B^*$, it is not necessarily convex. The outer normal transform of $B$ can thus be identified with the convex hull of $\widehat B$, and so, when $\widehat B$ is not convex, will possess a different combinatorial structure than $B^*$.  

Indeed, counterexamples to the problem of inscribing convex polytopes of a given combinatorial type in spheres, \cite{PZ}, are of this form.  For example, the dual to the truncated tetrahedron, known as the triakis tetrahedron, is not inscribable in a sphere.  The flattened version of the spherical dual to a truncated tetrahedron is a cube with diagonals that bisect each square into a pair of triangular facets, and form the edges of an interior tetrahedron.  Both the spherical dual and the resulting flattened cube with diagonals have the same combinatorial type as the triakis tetrahedron.  However, the flattened cube, while inscribed in the unit sphere, is not a convex polyhedron since it has pairs of coplanar triangular facets possessing a common normal.  It is, of course, the set-theoretic boundary of a convex subset of $\R^3$, namely the inscribed solid cube, whose cubical boundary (without the diagonals) can be identified as the outer normal transform of the original truncated tetrahedron, and is not combinatorially equivalent to the triakis tetrahedron. Furthermore, slightly perturbing the original truncated tetrahedron leads to a perturbed spherical dual and a perturbed cube with diagonals that is inscribed in the sphere, again both having the same combinatorial type as the triakis tetrahedron.  However, although its triangular faces are no longer coplanar, the resulting polyhedron is not the boundary of a convex subset of $\R^3$, and hence not equal to its outer normal transform, which is the convex hull of this nonconvex perturbed cube.  All this is a necessary consequence of the non-inscribability of the triakis tetrahedron. 

In general, if the flattened spherical dual of a polytope is convex then it has to coincide with its outer normal transform, which is then, by the above remarks, combinatorially equivalent to the dual polytope.  On the other hand, if it is not convex then its convexification, which is the outer normal transform, cannot be combinatorially equivalent to the dual.  Thus, we have established the following result.

\begin{proposition} \label{spherical_dual} Let $B \subset \Rd$ be a convex polytope of dimension $d-1$.  Then the outer normal transform of $B$ is combinatorially equivalent to the dual polytope $B^*$ if and only if the flattened spherical dual of $B$ is convex.
\end{proposition}


Finally, for later purposes, we will introduce some useful open subsets of the normal spherical complex \eq V.  If $F\subset B$ is a facet with outwards unit normal $\n_F \in \Sdo$, so $S_F = \{\n_F\}$, set
\Eq{WF1}
$$W_F = S_F \>\vee\> \bigvee_{G \subsetneq F} \> S_G^\circ  $$
where the union is over the proper subfaces $G \subsetneq F$.  On the other hand, if $F \subset B$ is a face with $1 \leq  \dim F < d-1$, set
\Eq{WF2}
$$W_F = \bigvee_{G \subseteq F} \> S_G^\circ  .$$

\begin{lemma} \label{WF} Under the above definitions, $W_F$ is a relatively open subset of $\Sdo$. \end{lemma}

\Proof
This follows from the fact that the corresponding union of normal cones
$$V_F = \{0\} \> \vee \> \bigvee_{G \subseteq F} \> N_G^\circ  $$
 is an open cone and $W_F = V_F \capz \Sdo$.  Indeed, one can use a perturbed version of the Supporting Hyperplane Theorem that says that  a perturbed supporting hyperplane remains supporting at some point in its support.  In more detail, if $H$ is a supporting hyperplane such that $H\capz B = F$ where $F$ is a face, and $\widetilde H$ is any sufficiently small perturbation of $H$, then $\widetilde H \capz B = G$ for some subface $G \subseteq F$.  Keep in mind that the subface could be a vertex.
\qed  

%

\medskip

The last result of this section is a technical construction, that is key to our proof of the Main Theorem \ref{main_theorem}.  The reader may wish to skip it for now, and return once the proof is underway.


\begin{proposition} \label{Nstar}  
Let $F \subset B$ be a face of dimension \hbox{$1 \leq m \leq d-1$}. Let $S_F$ be its normal spherical polytope and $W_F\subset \Sdo$ the open subset given by Lemma \ref{WF}. Let $\subs Gk$ be its $(m-1)$-dimensional subfaces, so that $\partial F = \bigcup_{i=1}^k G_i$.  Similarly, let $S^\circ_i = S^\circ_{G_i}$  and $W_i = W_{G_i}$.  

Suppose $N \subset W_F \subset \Sdo$ be a connected $m$-dimensional submanifold such that either {\rm ({\sl a\/})} if $F$ is a facet, of dimension $d-1$, with unit outwards normal $\n_0 = \n_F$, then $N$ is an open neighborhood of $\n_0$, or {\rm ({\sl b\/})} if $1 \leq m = \dim F < d-1$, then $N$ intersects $S_F^\circ$ transversally at a single point $\n_0 \in N \capz S_F^\circ$.  
Then if $\widetilde N \subset N$ is a sufficiently small open contractible submanifold with $\n_0 \in \widetilde N$, which implies  $ \n_0 \in \partial (\widetilde N \capz S_i^\circ)$ for all $\rg ik$,  we can decompose its boundary $\partial \widetilde N = \bigcup_{i=1}^k L_i$ into the union of $(m - 1)$-dimensional submanifolds that only overlap on their boundaries, meaning $\barr L_i \capz \barr L_j = \partial L_i \capz \partial  L_j $ whenever $i \ne j$, with the property that each $L_i \subset W_i$ intersects $S^\circ_i$ transversally at  a single point $\n_i \in L_i \capz S^\circ_i = \partial \widetilde N \capz S^\circ_i$.
\end{proposition}

\Proof
Choose $r>0$ sufficiently small so that the relatively open submanifold 
$N_r = \set{\n \in N}{\norm{\n - \n_0} < r}$
has boundary
$\partial N_r = \set{\n \in N}{\norm{\n - \n_0} = r}$. 
Moreover, reducing $r$ if necessary, we claim that $\partial N_r$ intersects each $S^\circ_i$ transversally at a single point $\n_i \in  \partial N_r\capz S^\circ_i  $.  Indeed, in a small neighborhood $\n_0 \in U$ we can choose local coordinates centered at $\n_0$ such that, locally, $S_F^\circ \capz U$ is a $(d-m-1)$-dimensional subspace, $N_r \subset U$ is a transverse $m$-dimensional subspace, while $S^\circ_i\capz U$ is a $(d-m)$-dimensional half space with local boundary $\partial S_i\capz U = S_F^\circ\capz U$, from which the preceding claim is evident.

We now set $\widetilde N = N_r$. Since $S_F \subset \partial S_i$, this immediately implies $ \n_0 \in \partial (\widetilde N \capz S_i^\circ)$.
  The final task is to decompose $\partial \widetilde N= \bigcup_{i=1}^k L_i$ as in the statement of the Proposition. It is reasonably clear that there are many ways to do this, but for definiteness here is one possible construction.  First we note that, by \eqs{WF1}{WF2}, either
$$\qeq{W_F \setminus \{\n_F\} = \bigcup_{i=1}^k\ W_i,\rox{or} \ W_F \setminus S_F^\circ = \bigcup_{i=1}^k\ W_i,}$$
according to whether $F$ is a facet or not.  We thus, for each $\rg ik$, need to choose $L_i \subset \widetilde N \capz W_i$ with the requisite properties.

First, define the closed subset $\barr L_i \subset \partial \widetilde N$ to be the set of all $\n \in \widetilde N \capz W_i$ such that if $\n \in \widetilde N \capz S_H^\circ$ for some adjacent subface $H \subsetneq G_i$,  then $\dist(\n,\n_i) \leq \dist(\n,\n_j)$ for all other adjacent $(m-1)$-dimensional subfaces $G_j$, meaning that $H \subsetneq G_j$. Clearly $N = \bigcup \barr L_i$ and, moreover, $\barr L_i$ and $\barr L_j$ only overlap on their common boundary, which could be either part of a boundary of an $S_H$ or a point $\n \in N \capz S_H^\circ$ that is equidistant to $\n_i$ and $\n_j$.  We then set $L_i = \barr L_i^\circ$ to be its interior\footnote{It may happen that the closure of $L_i$ is strictly contained in $\barr L_i$; this can occur if there exist $\n_l$ associated with nonadjacent faces $G_l$ that lie closer to the points $\n \in \widetilde N \capz S_H^\circ$ than those in any adjacent face $G_i$.  But this does not affect the construction since every point in $\barr L_i \setminus L_i$ is contained in the boundary of some $L_j$.} relative to $ \partial \widetilde N$.  Since $\n_i \in L_i^\circ$, transversality of $\partial \widetilde N$ to $S_i^\circ$ at $\n_i$ immediately implies the same for the relatively open submanifold $L_i$.  We conclude that the resulting submanifolds  satisfy the required conditions. \qed

\section{The Gauss Map of a Convex Polytope} \label{Gauss}

We claim that the preceding construction can be identified with a form of the inverse of the Gauss map of the boundary of the convex hull $B = \partial K$. Recall, \cite{GAS}, that the \is{Gauss map} of a smooth closed hypersurface, i.e., a $(d-1)$-dimensional oriented submanifold $M\subset \Rd$, is
\Eq{Gauss}
$$\qeq{\gamma_M \colon M \longrightarrowx \Sdo,\\\gamma (\y) = \n_\y,\\\y \in M,}$$
where $\n_\y$ denotes the unit outward normal to $M$ at $\y$. If $M$ is convex, then its Gauss map is one-to-one and onto, with smooth inverse $\gamma_M^{-1} \colon \Sdo \to M$.  

We are interested in the convergence of the Gauss maps $\gamma _{M_{\varepsilon}}$ associated with a parametrized family of smooth closed convex hypersurfaces $M_{\varepsilon}$, for $\varepsilon >0$, that converge to the piecewise linear convex hypersurface (polytope) $B = \partial K$ as $\varepsilon \to 0^+$.  
Convergence of the Gauss maps will be in the sense of  set-valued functions, as we now describe. 

In general, a \is{set-valued function}, also known as \is{multi-valued functions}, from a space $D$ to a space $Y$ means a mapping $F$ from $D$ to the power set $2^Y$, i.e., the set of subsets of $Y$, \cite{AF,AO}. In other words, the image of $x \in D$ is a subset  $F(x) \subset Y$.   More generally, a set-valued function maps subsets of its domain to subsets of its range in the evident manner. We say that $F$ \emph{has closed values} if $F(x)$ is a closed subset of $Y$ for all $x \in D$. In particular, any ordinary function $f \colon D \to Y$ can be viewed as a set-valued function, with closed values, by identifying the image $y = f(x)$ of a point $x \in D$ with the singleton set $\{y \} \subset Y$. The \emph{range} $R \subset Y$ of $F$ is the union of all the images of points in its domain $D$, so $R = F(D)$.  

Here is a simple example of convergence of set-valued functions.

\medskip

\begin{example}\label{arctan}  
\rm 
Consider the ordinary functions
\Eq{arctan}
$$f_{\varepsilon}(x) = \frac 2\pi \,(1- \varepsilon ) \arctan \frac x{\varepsilon} \quad \rox{for} x \in \R, \ \ \varepsilon > 0.$$
In the usual function-theoretic sense of convergence,
$$\lim_{\varepsilon \to 0^+} f_{\varepsilon}(x) = \sign x = \mcases{1,& x > 0,\\0,&x = 0,\\-1,&x <0.}$$
Thus, for almost every point $x \in \R$, the value of $f_{\varepsilon}(x)$ converges to either $-1$ or $1$.
However, if you look at their graphs as subsets of $\R^2$, they converge, as sets, to the curve consisting of the union of the three line segments 
$$\set{(x,-1)}{x \leq 0} \cupz \set{(0,y)}{-1 \leq y \leq 1} \cupz \set{(x,1)}{x \geq 0} .$$ 
We can interpret this curve as the graph of the set-valued function 
\Eq{setF}
$$\widehat f \colon \R \longrightarrowx 2^{\R} \roq{given by} \widehat f(x) = \xcases{\{1\},& x > 0,\\[-1,1],&x = 0,\\ \{-1\},&x <0.}$$
The domain of $\widehat f$ is $D = \R$ and its range is the interval $R = [-1,1] = \widehat f(\R)$. Furthermore, the ranges $R_\varepsilon = f_\varepsilon (\R) = (-1 + \varepsilon , 1 - \varepsilon )$ of the functions \eq{arctan} are open intervals that converge, as sets, to the closed interval $[-1,1]$ forming the range of the limiting set-valued function.
\end{example}

\medskip


In general, given spaces $D,R$, let $ D \times R$ denote their Cartesian product, and $\pi_D \colon  D \times R \to D$ and $\pi_R \colon  D \times R \to R$ the standard projections.  Any subset $S \subset D \times R$ which projects onto both 
$D = \pi_D(S)$ and $R = \pi_R(S)$ defines a set-valued mapping $F$ with domain $D$ and range $R$, given by $F(x) = \pi_R( S \cap \pi_D^{-1}\{x\})$.   Its inverse $F^{-1}$ is also a set-valued mapping from $R$ to $D$, given by $F^{-1}(y) = \pi_D( S \cap \pi_R^{-1}\{y\})$, \cite{AF}.  
For the above example \eq{setF}, $\widehat f^{-1}\colon [-1,1] \to 2^{\R}$ is given by
\Eq{setFi}
$$\widehat f^{-1}(y) = \xcases{(-\infty ,0],& y=-1,\\\{0\},&-1< y < 1,\\[0,\infty ),&y = 1.}$$
Note that any ordinary function thus has a set-valued inverse. 

In convex analysis, the normal cone \eq{NS} is often viewed as a set-valued mapping, that maps a point $\x \in B = \partial K$ to its normal cone $N_\x$. (One can, of course, extend it to all of $K$ but the values on the interior $K^\circ$ are trivial.) Here we consider our normal spherical polytope construction as a set-valued mapping $\gamma _B$ from $B$ to $\Sdo$, mapping a point $\x \in B$ to its normal spherical polytope: $\gamma _B(\x) = S_\x \subset \Sdo$. 

Suppose that $M_{\varepsilon}$ for $\varepsilon > 0$ are a parametrized family of smooth closed convex hypersurfaces converging uniformly to the boundary of the convex hull: $M_{\varepsilon} \to B$ as 
$\varepsilon \to 0^+$.  Then their Gauss maps $\gamma _{M_{\varepsilon}} \to \gamma _B$ converge, in the sense of set-valued functions, to the set-valued normal spherical polytope map.  We can thus identify the normal spherical polytope map constructed above as the \emph{Gauss map of a convex polytope, i.e., a piecewise linear convex hypersurface.}

As for our construction, the Main Theorem \ref{main_theorem} shows that the functions $\ffe \colon \Sdo \to \Rd$ converge, in the set-valued sense, to the inverse of the Gauss map $\gamma _B^{-1}\colon \Sdo \to B$ associated with the boundary of the convex hull.  On the other hand, the $\ffe $ are certainly not inverse Gauss maps themselves. Moreover, simple examples, e.g., that in Figure \ref{fig3}, show that the image $\ffe (\Sdo)$ is not in general a convex hypersurface. On the other hand, it might be worth investigating the set-theoretic convergence of their possibly multi-valued Gauss maps.

Finally, we note that the concept of continuity does not extend straightforwardly to set-valued functions, \cite{AF}. The most important analog is contained in the following definition.

\begin{definition}\label{hemicontinuous} Let  $D,R$ be topological spaces. A set-valued mapping $f \colon D \to R$ is called  \emph{upper hemicontinuous at $x_0 \in D$} if and only if, for any open neighborhood $V$ of the set $f(x_0)$, there exists a 
neighborhood  $U$ of $x_0$ such that $f(x) \subset V$ for all $x \in U$. We say that $f$ is \emph{upper hemicontinuous} if it is upper hemicontinuous at every $x_0 \in D$.  \end{definition}

It is straightforward to verify that Example \ref{arctan} satisfies the upper hemicontinuity condition.

\smallskip

%
%

\emph{Warning\/}: A few authors, including \cite{AF}, use the expression ``upper semicontinuous'' instead of ``upper hemicontinuous''.  However the latter terminology seems to be more accepted by the broader community, particularly as it is not in conflict with the notion of semicontinuity of ordinary functions.


\section{Some Computational Lemmas} \label{lemmas}

Before launching into the proof of the Main Theorem \ref{main_theorem}, let us collect together some elementary computational lemmas for the functions used to form the maps $\ffe$ defined in \eq{fe}.

Recalling \eq{cij} and \eq{ci}, let us set
\Eq{cij0}
$$\qeq{c_{ij}(\n) = \lim_{\varepsilon \to 0^+} \, c_{ij}(\varepsilon,\n) = \max\{ 0, -\ip{\n}{\n_{ij}}\},\\\n \in \Sdo.}$$
and
\Eq{ci0}
$$\qeq{c_i(\n) = \lim_{\varepsilon \to 0^+} \, c_{i}(\varepsilon,\n) = \prod_{\substack{1 \leq j \leq n \\ j \neq i}} c_{ij}(\n),\\\n \in \Sdo.}$$
We can thus write
\Eq{bi}
$$c_i(\varepsilon, \n) = c_i(\n) + \varepsilon \,b_i(\varepsilon, \n),$$
where $b_i$ is a polynomial in $\varepsilon $ of degree $n-2$.  
In view of \eq{Sio} and \eq{SiS}, $c_i(\n) > 0$ if and only if $\n \in S^\circ_i$. Thus, $c_i(\n) = 0 $ for all $\rg in$ if and only if 
$\n \in \Q= \Sdo  \setminus S$, where $S = \bigvee_{i=1}^\nv \>S^\circ_i $
is the disjoint union of the interiors of the normal spherical polytopes associated with the vertices $\subs \x \nv$ of $B$.  
We have thus established the following result.

\begin{lemma} \label{main_lemma} Given $\n \in \Sdo$, either all $c_i(\n) = 0$, or precisely one $c_i(\n) > 0$, and the rest are all zero. Moreover, in the latter case, $\x_i$ is a vertex. \end{lemma}

Indeed, if $\n \in S^\circ _i$, then, referring to \eqs{ci}{bi}, $c_i(\n) > 0$ and $b_i(\varepsilon,\n) > 0$ for all $\varepsilon > 0$, whereas
\Eq{cjoe}
$$c_j(\varepsilon,\n) = \varepsilon^{1+k_j} \, a_j(\varepsilon,\n) \rox{for} j \ne i, $$
with $a_j(\varepsilon,\n) > 0$ for  $\varepsilon > 0$.  Here the nonnegative integer $k_j$ denotes the number of points $\x_k$ with $k \ne i,j$ that satisfy the distance inequality 
$\dist(\x_k, P_i(\n)) \leq \dist(\x_j, P_i(\n))$ where $P_i(\n) = \x_i + \n^\perp$, with $\n^\perp = \set{\y \in\Rd}{\ip \y \n = 0}$ denoting the affine hyperplane orthogonal to $\n$ passing through $\x_i$.  Note that there is either one or no value of $j$ for which $k_j = 0$.


Let us finish this section by establishing a more detailed version of Lemma \ref{main_lemma}, valid for an arbitrary face $F \subset B$. 

\begin{lemma} \label{face_lemma} 
Let $F \subset B$ be a face of dimension $1 \leq m \leq d-1$, with vertices $\subs \x k$.  Suppose $F$ contains $l$ additional {\rm (}non-vertex\/{\rm)} points $\x_{k+1},\ldotsx, \x_{k+l} \in F^\circ$, where $l$ may be zero.  Then, given $\n \in S^\circ_F$, 
\Eq{cie}
$$\eeq{c_i(\varepsilon,\n) = \varepsilon ^{k+l-1} d_i(\n) + \varepsilon ^{k+l}p_i(\varepsilon ,\n),& \rg i{k+l},\\
c_j(\varepsilon,\n) = \varepsilon ^{k+l}p_j(\varepsilon ,\n), & \rgs j {k+l+1}n,}$$
where $d_i(\n) > 0$, while $p_1(\varepsilon ,\n), \ldots, p_n(\varepsilon ,\n)$ are polynomials in $\varepsilon $.
\end{lemma}

\Proof
Let $\n \in S^\circ_F$.  The proof follows from the fact that $c_{ij}(\varepsilon ,\n) = \varepsilon $ whenever $\x_i,\x_j \in F$, so that $\ip \n{\n_{ij}} =  0$.  On the other hand, $\ip \n{\n_{ij}} <  0$ whenever $\x_i \in F$ and $\x_j \not\in F$, which, vice versa, implies $c_{ji}(\varepsilon ,\n) = \varepsilon $. The proof is completed by recalling the definition \eq{ci} of $c_i(\varepsilon,\n)$. \qed

\section{Proof of the Main Theorem} \label{main}

Now we turn to the proof of the Main Theorem \ref{main_theorem}, on the convergence, as $\varepsilon \to 0^+$, of the hypersurfaces $\ffe(\Sdo)\subset K^\circ \subset \Rd$ to the boundary of the convex hull $B = \partial K$ of the point configuration $\X$.


If we formally set $\varepsilon =0$ in the preceding definition \eq{fe} of the map $\ffe $,  Lemma \ref{main_lemma} implies
\Eq{limf}
$$\Lim \varepsilon {0^+} \ffe(\n) = \mcases{\x_i,& \n \in S^\circ_i,\\\roh{undefined},& \textstyle \n \in \Q = \Sdo  \setminus \bigvee_{i=1}^\nv \>S^\circ_i .}$$
Thus, for almost every point $\n \in \Sdo$, the images $\ffe(\n)$ converge to one of the vertices of the convex hull.  However, as in  Example \ref{arctan}, \emph{this does not imply that, as a set, $\ffe (\Sdo)$ converges to the set of vertices $V = \bsubs \x \nv$}.  Our goal is to prove that the images $\ffe (\Sdo)$ converge, as sets, to the entire boundary $B$ as $\varepsilon \to 0^+$.  Specifically, we will show:

\bitems{Given any neighborhood $W \supset B$, no matter how small, we can find $\varepsilon_0 > 0$ such that $\ffe (\Sdo) \subset W$ for any $0 < \varepsilon < \varepsilon_0$.\\ 
Given any $\x \in B$, there exist points $\y_\varepsilon \in \ffe (\Sdo)$ for $\varepsilon > 0$ such that $\y_\varepsilon \to \x$ as $\varepsilon  \to 0^+$.  
}

In the language of set-theoretic limits, \cite{AF}, the first statement shows that the \emph{outer limit} of the sets $\ffe (\Sdo)$ is a subset of $B$. The second statement proves that $B$ is a subset of their \emph{inner limit}. Since the inner limit is always a subset of the outer limit, this then implies that the inner limit and outer limit coincide and are equal to $B$.

First, recall that, for $r > 0$, the \emph{$r$-neighborhood} $U_r$ of a subset $D \subset \Rd$ is the set of points that are a distance less than $r$ (in the Euclidean norm) from $D$, i.e., $U_r = \set{\x \in \Rd}{\dist(x,D) < r}$.  In what follows, when we refer to an \emph{$\O(\varepsilon )$ neighborhood} of a set, by which we mean an $\varepsilon $ dependent system of $r$-neighborhoods in which, for $\varepsilon $ sufficiently small, $r = c \, \varepsilon $ for some unspecified constant $c$.

In order to understand our set-theoretic limit, we will investigate the behavior of the images $\ffe(A)$ of certain subsets $A \subset \Sdo$, gradually building up to the entire sphere.  
Let us begin with the simplest case: the images of a curve $C \subset \Sdo$.  If $C \subset S^\circ_i$ is entirely contained in the interior of the normal spherical polytope associated with a vertex $\x_i$ for some $1 \leq i\leq k$, then, by \eq{limf}, $ \ffe(C) \to \{\x_i\}$ as $\varepsilon \to 0^+$.

The next simplest case is when the curve $C$ is contained in the union of two adjacent vertex spherical polytopes. Thus, by relabeling, let $\x_1,\x_2$ be adjacent vertices of $B$.  Let 
$$E =\set{\lambda_1  \,\x_1 + \lambda _2\,\x_2 }{\lambda _1, \lambda _2 \geq  0, \ \ \lambda _1 + \lambda _2 = 1} \subset B$$
denote the edge connecting $\x_1$ to $\x_2$.  
Suppose that its interior contains $l\geq 0$ additional points in our configuration, which we number as $\x_3,\ldots, \x_{l+2} \in E^\circ$, while the remaining points $\x_{l+3},\ldots, \x_n \in K \setminus E$.  We note that we can also write, redundantly,
\Eq{E12}
$$E =\set{\Sum i{l+2} \lambda_i  \,\x_i }{\lambda _i \geq  0, \ \ \Sum i{l+2} \lambda_i = 1} .$$

Let $S_1,S_2 \subset \Sdo$ be the normal spherical polytopes associated with $\x_1,\x_2$, respectively, while $S_E = \partial S_1 \capz \partial S_2$ is the normal spherical polytope associated with the edge $E$. Thus $S^\circ_1 , S^\circ_2$ are open subsets of $\Sdo$, while $S^\circ_E$ is a $(d-2)$-dimensional submanifold.  Consider a curve $C \subset S^\circ_1 \cupz S^\circ_2 \cupz S^\circ_E \subset \Sdo$ such that one endpoint of $C$ lies in $S^\circ_1$ and the other lies in $S^\circ_2$, which, by connectivity, imply $C\capz S_E^\circ \ne \emptyset$.
Our goal is to prove that the image curves $\ffe(C)$ converge, as sets, to the edge $E$. 

Now, if $\n \in C \capz S_1^\circ $, Lemma \ref{main_lemma} combined with equations \eq{li} and \eq{cjoe} imply
\Eq{l1j}
$$\xeq{\lambda _1(\varepsilon,\n) = 1 + \varepsilon \,q_1(\varepsilon ,\n),&\lambda _j(\varepsilon,\n) = \varepsilon \,q_j(\varepsilon ,\n), & j=2,3,\ldots, n,}$$
where $\subs qn$ are rational functions of $\varepsilon $ depending continuously on $\n \in C$. Thus, \eq{l1j} re-establishes the fact that all of the points in $\ffe(C \cap S^\circ_1)$ converge to the vertex $\x_1$  as $\varepsilon \to 0^+$. 
A similar statement holds for $\n \in C\capz S_2^\circ$:
\Eq{l2j}
$$\xeq{\lambda _2(\varepsilon,\n) = 1 + \varepsilon \,q_2(\varepsilon ,\n),&\lambda _j(\varepsilon,\n) = \varepsilon \,q_j(\varepsilon ,\n), & j=1,3,\ldots, n.}$$
  Finally, if $\n \in C \capz S^\circ_E$, in view of \eqss{li}{Delta}{cie}, we have
\Eq{liq}
$$\lambda _i(\varepsilon ,\n)= \mcases{\frac{d_i(\n) + \varepsilon \,p_i(\varepsilon,\n)}{D(\n)+ \varepsilon \,P(\varepsilon,\n)},&\rg i{l+2},\\\frac{\varepsilon \,p_i(\varepsilon,\n)}{D(\n)+ \varepsilon \,P(\varepsilon,\n)},&\rgs i{l+3}n,} $$
where
$$\qeq{D(\n) = \Sum i{l+2}d_i(\n) > 0,\\P(\varepsilon,\n) = \Sum i{n}p_i(\varepsilon,\n).}$$
Comparing with \eqss{l1j}{l2j}{liq}, we find that for any $\n \in C$,
\Eq{felx}
$$\qeq{\Sum i{l+2}\lambda _i(\varepsilon,\n) = 1 + \varepsilon \, Q(\varepsilon,\n),\\
\ffe (\n) = \Sum i{l+2}\lambda _i(\varepsilon,\n) \x_i + \varepsilon \, R(\varepsilon,\n),}$$
where both $Q$ and $R$ are continuous functions of $\n \in C$, including when $\n \in S_E^\circ$, and rational functions of $\varepsilon $ with nonvanishing denominator.  Since $C \subset \Sdo$ is compact, they can thus be bounded by an overall constant independent of $\varepsilon \in (0, \varepsilon _0]$. This holds even at the singular point $\n \in C \capz S_E^\circ$ when there is cancellation of powers of $\varepsilon $ in numerator and denominator, whence \eq{liq}. Thus, comparing with \eq{E12}, we deduce that, for $0  < \varepsilon \leq \varepsilon _0$, the images $\ffe(C)$ lie in an $\O(\varepsilon )$ neighborhood $\Ue$ of the edge $E$. This immediately implies that the limiting set is contained within the edge: $\lim_{\varepsilon \to 0^+} \ffe(C) \subset E$.  
The remaining task is to prove that every point in $E$ is contained in the limit, and therefore $\lim_{\varepsilon \to 0^+} \ffe(C) = E$.



We already know that both endpoints $\x_1,\x_2$ are contained in the limiting set.  Thus, our remaining task is, given a point $\x \in E^\circ$, to find points $\y_\varepsilon \in \ffe(C)$ that converge to $\x = \lim_{\varepsilon \to 0^+} \y_\varepsilon$.  Although it is possible to do this by a careful analysis of the underlying formulae, we prefer, for later purposes, to use a simple topological proof. 

To this end, let $Z_\x$ be the affine hyperplane passing through $\x$ that is orthogonal to $E$, and define $Z_{\x,\varepsilon} = Z_\x \cap \Ue$.  We claim that there exists $\y_ \varepsilon \in \ffe(C) \capz Z_{\x,\varepsilon} $. If true, then we have produced the desired points. To prove the claim, observe that $\Ue \setminus Z_{\x,\varepsilon}$ consists of two disjoint open subsets, say $U^1_{\x,\varepsilon},U^2_{\x,\varepsilon}$ with $\x_i \in U^i_{\x,\varepsilon}$ for $i=1,2$.  
Moreover, since we know that all the points in $\ffe(C\capz S^\circ_i)$ converge to $\x_i$, if we choose $\varepsilon$ sufficiently small, then $\ffe(C)  \capz U^i_{\x,\varepsilon} \ne \emptyset$.  Therefore, $\ffe(C) \capz Z_{\x,\varepsilon} = \emptyset$  would contradict the connectedness of $\ffe(C)$. This contradiction establishes the above claim.  
We thus conclude that, as sets
\Eq{setconv}
$$\ffe(C) \longrightarrow E \rox{as} \varepsilon \longrightarrow 0.$$

For later purposes, we need slightly more than mere set-theoretic convergence \eq{setconv}.  Namely, we require the existence of a continuous set-valued homotopy that connects the images of $\ffe\colon C \to K$ for $\varepsilon > 0$ to a set-valued map $\widehat \ff_0\colon C \to 2^K$ with range equal to the edge $E = \widehat \ff_0(C)$, a model being Example \ref{arctan}. Rather than write down an explicit formula for this homotopy, we will instead construct its graph.  

Consider the graph
$$\Gamma = \set{\bpa{\varepsilon , \n, \ffe(\n)}}{0 < \varepsilon \leq \varepsilon _0, \  \n \in C} \subset (0,\varepsilon _0\,]\times C \times K$$
of the map $F(\varepsilon , \n) = \ffe(\n)$ for $0 < \varepsilon \leq \varepsilon _0$ and $\n \in C$.   
Let $\barr \Gamma = \roh{Clos}\> \Gamma$ be its closure in $[0, \varepsilon _0] \times C \times K$.  According to the preceding proof, $\barr \Gamma$ is the graph of the set-valued map $F \colon [0, \varepsilon _0] \times C \to 2^K$ given by
\Eq{FE}
$$F(\varepsilon ,\n) = \mcases{\ffe(\n),& \varepsilon > 0,\\\x_1,&\varepsilon = 0,\ \ \ \n \in C \capz S_1^\circ,\\\x_2,&\varepsilon = 0,\ \ \ \n \in C \capz S_2^\circ,\\E,&\varepsilon = 0,\ \ \ \n \in C \capz S_E^\circ,}$$
its final value being the entire edge $E\subset K$.  
Moreover, since $\barr \Gamma$ is closed and $K$ is compact Hausdorff, the Closed Graph Theorem for set-valued functions, \cite[Prop.~1.4.8]{AF}, implies that the set-valued function $F$ is upper hemicontinuous, as per Definition \ref{hemicontinuous}. Thus, for all $0 < \varepsilon \leq \varepsilon _0$, \eq{FE} defines an upper hemicontinuous homotopy  from each $\ffe \colon C \to K$ to the set-valued map $\widehat \ff_0\colon C \to 2^K$ with $\widehat \ff_0(\n) = F(0,\n)$, whose range $\widehat \ff_0(C)$ is the edge $E$.

\medskip

\emph{Remark\/}: An alternative approach, that avoids set-valued homotopies and, later, set-valued homology, is to ``tilt'' the subset $\barr\Gamma $ so that it becomes a graph by introducing new coordinates on the Cartesian product space $[0, \varepsilon _0] \times C \times K$.  However, this is more technically tricky to accomplish in the higher dimensional cases to be handled below, and the set-theoretic approach provides a cleaner path to the proof.

\medskip

The remainder of the proof works by induction on the dimension of the face $F$.  Thus, the next case is that of a two-dimensional face $F \subset B \subset \R^d$.  The main steps of the proof in this situation will then be straightforwardly adapted to any higher dimensional face. Let $\subs \x k $ be the vertices of $F$ and let $\subs Ek$ be its edges.   We label the vertices and edges so that $E_j$ connects $\x_j$ to $\x_{j+1}$, with indices taken modulo $k$ throughout, whence $\x_{k+1} = \x_1$.  Thus $E = \bigcup_{j=1}^k E_j = \partial F$ is the polygonal boundary of $F$.  We assume that there are $l\geq 0$ additional points $\x_{k+1},\ldots, \x_{k+l} \in F\setminus \bsubs \x k$, while the remaining points in the configuration $\x_{k+l+1},\ldots, \x_n \in K \setminus F$.  Keep in mind that $F$ is convex.

Let $S_i, \widetilde S_j,S_F$ be the normal spherical polytopes of $\x_i, E_j, F$, respectively, so that $\widetilde S_j \subset \partial S_j \capz \partial S_{j+1}$ and $S_F \subset \partial  \widetilde S_j$ for all $\rg jk$. 
Let $W_F \subset \Sm^{d-1}$ be the open set \eqs{WF1}{WF2}, and let $\widetilde N \subset N \subset W_F$ be the two-dimensional submanifolds satisfying the hypotheses of Proposition \ref{Nstar}.
As in \eq{felx}, applying Lemma \ref{face_lemma}, we find
\Eq{Felx}
$$\qeq{\Sum i{k+l}\lambda _i(\varepsilon,\n) = 1 + \varepsilon \, Q(\varepsilon,\n),\\
\ffe (\n) = \Sum i{k+l}\lambda _i(\varepsilon,\n)\, \x_i + \varepsilon  R(\varepsilon,\n),}$$
where both $Q$ and $R$ are continuous functions of $\n \in N$, and rational functions of $\varepsilon $ with nonvanishing denominator. These formulae again imply that the images $\ffe(N)$ lie in an $\O(\varepsilon )$ neighborhood $\Ue$ of the face $F$, and hence $\lim_{\varepsilon \to 0^+} \ffe(N) \subset F$. The remaining task is to prove that every point in $F$ is contained in the limit, a result that requires a more sophisticated topological argument than in the curve case. 

For this purpose, we replace $N$ by $\widetilde N$. Clearly, if we can prove $\lim_{\varepsilon \to 0^+} \ffe(\widetilde N) = F$, by the preceding result the same is true of $N \supset \widetilde N$. According to Proposition \ref{Nstar},  $\widetilde N \capz S_j^\circ \ne \emptyset$ and $\widetilde N \capz \widetilde S_j^\circ \ne \emptyset$ for all $\rg jk$ and either $\n_F \in \widetilde N$ when $d=3$, where $\n_F$ is the unit outward normal to the polyhedral facet $F$, or $\widetilde N \capz S_F^\circ \ne \emptyset$ when $d > 3$.  Moreover, the boundary $L = \partial \widetilde N$ can be decomposed into nonoverlapping curves $\subs Lk$ that satisfy $L_j \subset S_j^\circ \cupz \widetilde S_j^\circ \cupz S_{j+1}^\circ $, again modulo $k$.  Let $\bco{\n_j} =  L_{j-1} \capz L_j \subset S_j^\circ$ denote the common endpoints of adjacent curves in $L = \partial \widetilde N$.

Let us set $\I = [\,0,\varepsilon _0\,]$ for $\varepsilon _0 > 0$ sufficiently small. 
According to the preceding curve proof, $\ffe(L_j) \to E_j$ as sets and, moreover, there exists an upper hemicontinuous homotopy (of set-valued mappings) from each $\ffe\colon L_j \to K$ for all $0 < \varepsilon \leq \varepsilon _0$ to the set-valued limit $\widehat \ff_0 \colon L_j \to 2^K$ with range equal to the edge $E_j = \widehat \ff_0(L_j)$. The graph of this homotopy,
$$\Gamma _j \subset  \I \times L_j \times K \subset \I \times \Sm^2 \times  K,$$
is a closed subset of the indicated Cartesian product space.  

We now piece together these homotopy graphs to define a homotopy from $\ffe(L)$ to $E = \partial F$ whose graph is
\Eq{Gamma}
$$\Gamma = \bigcup_{j=1}^k \ \Gamma _j \subset \I \times L \times K \subset  \I \times \Sm^2 \times K.$$
Note that $\Gamma $ is a closed subset that defines the graph of an upper hemicontinuous function because each $\Gamma _j$ is closed and, moreover,
$$\Gamma _{j-1} \cap \bpa{\I \times \bco{\n_j} \times K} = \Gamma _j \cap \bpa{\I\times \bco{\n_j} \times K},$$
including when $\varepsilon = 0$, since $\ffe(\n_j) \to \x_j$, and so 
$$\Gamma _{j-1} \cap \bpa{\bco0 \times \bco{\n_j} \times K} = \bco{(0,\n_j,\x_j)} = \Gamma _j \cap \bpa{\bco0 \times \bco{\n_j} \times K}.$$

Given $\x \in F^\circ$, we seek $\y_\varepsilon \in \ffe(\widetilde N)$ that converge to $\x$ as $\varepsilon \to 0^+$. Let $Z_\x$ be the affine subspace of dimension $d-2$ passing through $\x$ that is orthogonal to $F$.  Define $Z_{\x,\varepsilon} = Z_\x \cap \Ue$. Again, if we can prove there exists $\y_\varepsilon \in  Z_{\x,\varepsilon} \capz \ffe(\widetilde N)$, we are done.  Suppose not, i.e., suppose that $\ffe(\widetilde N) \subset \Ue \setminus Z_{\x,\varepsilon}$. The idea is to demonstrate that this is topologically impossible due to the contractibility of $\widetilde N$, and hence contractibility of $\ffe(\widetilde N)$, whereas $\ffe(L) = \ffe(\partial \widetilde N)$, for $\varepsilon $ sufficiently small, defines a nontrivial homology class in $\Ue \setminus Z_{\x,\varepsilon}$.  

If we were dealing with ordinary mappings, this topological argument would be straightforward.  But because $\widehat \ff_0$ is a set-valued mapping, we will need some more sophisticated tools from set-valued algebraic topology to establish the contradiction.  We summarize the basic theory, based on a paper of Yongxin Li, \cite{Li}, in Appendix \ref{appB}.  In accordance with the notation introduced there, we use roman $H_n(X)$ to denote the standard $n$-th order singular homology groups of a topological space $X$, and calligraphic $\mathcal{H}_n(X, \mathcal{U})$ to denote the corresponding $n$-th order set-valued homology groups relative to a chosen open cover $\mathcal{U}$.  (As noted in the appendix, if one does not choose this cover carefully, the set-valued homology groups are all trivial, and would hence be useless for our purposes.)

In our situation, we select the particular open covering $\mathcal{V}$ of $\Ue \setminus Z_{\x;\varepsilon}$ consisting of all open sets of the form
\Eq{Vo}
$$ V = H \cap \bpa{\Ue \setminus Z_{\x;\varepsilon}} \rox{such that} \x \in \partial V,$$
where $H$ is an open half-space in $\R^d$. We claim that $\mathcal{V}$ satisfies Li's contractible finite intersection property, 
because the intersection of any finite collection of such open sets, if non-empty, is homeomorphic to the Cartesian product of an open $(d-2)$-dimensional ball with an open circular sector, meaning the intersection of an acute-angled open planar cone with the unit disk (a pizza slice), which is clearly contractible. 


\begin{figure} 
\centering
\begin{tikzpicture}[node distance=3.5cm, auto]
  \node (G) {};
  \node (O) [below of=G]{$H_1(\partial \widetilde N)$};
  \node (A) [right of=G]{$\mathcal{H}_1(\Ue \setminus Z_{\x;\varepsilon}, \,\mathcal{V})$};
  \node (B) [right of=O] {$\mathcal{H}_1(\Ue \setminus Z_{\x;\varepsilon}, \,\mathcal{V})$};
  \node (C) [below of=B]{$H_1(\Ue \setminus Z_{\x;\varepsilon})$};
  \node (E) [below of=O]{$H_1(\widetilde N) = \mathbf{0}$};
  \draw[->] (O) to node {$(\widehat \ff_{0})_*$} (A);
  \draw[->] (O) to node {$(\ff_{\varepsilon})_*$} (B);
  \draw[->] (O) to node [swap] {$\ff_{\varepsilon,*}$} (C);
  \draw[->] (A) to node {$\Id$} (B);
  \draw[->] (C) to node [swap] {$i_\sharp$} (B);
  \draw[->] (O) to node [swap] {$\iota_*$} (E);
  \draw[->] (E) to node {$\ff_{\varepsilon,*}$} (C);
\end{tikzpicture}
\caption{Commutative Diagram} \label{diagram} 
\end{figure}

The limiting set-valued function $\widehat \ff_0\colon L \to 2^{E}$, whose range is the polygonal boundary of the face $\widehat \ff_0(L) = E = \partial F$, is compatible with the open covering $\mathcal{V}$, because $\widehat \ff_0(\n)$ is either a vertex or an edge $E_j$. Moreover, when $\varepsilon >0$, the map $\ffe$ is continuous and single-valued, which implies trivially that its restriction to $L$ is compatible with any open covering of $\Ue \setminus Z_{\x;\varepsilon}$.

The family of maps $\{\ffe,\widehat \ff_0\}$ thus defines, by varying $\varepsilon$, 
an upper hemicontinuous homotopy of multi-valued functions with closed values. It follows from \cite[Prop.~6]{Li} that the upper triangle in Figure \ref{diagram} commutes.  
In the same figure, the square is divided into two triangles. It follows from the definitions that the top right triangle in the square commutes.
As noted in the Appendix \ref{appB}, the map $i_{\sharp}$ on the bottom right is an isomorphism. Finally, the map $\iota \colon L = \partial \widetilde N \to \widetilde N$ denotes the inclusion map, and so it is a standard fact from ordinary singular homology theory that the bottom left triangle commutes.

Let $0 \ne c = [\partial \widetilde N] \in H_1(\partial \widetilde N)$ be the homology class representing $\partial \widetilde N$, which is, in fact, a generator. We claim that\footnote{As in Appendix \ref{appB}, the parentheses indicate the induced maps on set-theoretic homology.} $(\widehat \ff_0)_*(c)$ is a non-zero element of $H_1(\Ue \setminus Z_{\x;\varepsilon}, \,\mathcal{V})$. Indeed, $(\widehat \ff_0)_*(c) = i_\sharp\bpa{[E]} \ne 0$ since the homology class $[E] = [\partial F] \in H_1(\Ue \setminus Z_{\x;\varepsilon})$ is nonzero and $i_\sharp$ is an isomorphism.
This thus proves the claim that  
$$(\ffe)_*(c) = (\widehat \ff_0)_*(c) \neq 0, \rox{and hence}\ff_{\varepsilon ,*}(c) =  i_{\sharp}^{-1}\bbk{(\ffe)_*(c)} \neq 0. $$
On the other hand, the bottom left triangle in the square in figure \ref{diagram} shows that $\ff_{\varepsilon ,*}$ vanishes identically on $H_1(\partial \widetilde N)$, so that $\ff_{\varepsilon ,*}(c) = 0$, thus leading to the desired contradiction and thus establishing the existence of $\y_\varepsilon \in \ffe(\widetilde N)$.
This finishes the proof that every point of $F$ belongs to the inner limit of $\ffe(\widetilde N)$, as $\varepsilon \to 0$. We conclude that both $\ffe(\widetilde N)$ and $\ffe(N) \to F$ as sets as $\varepsilon \to 0^+$.

Finally, to establish the existence of a set-valued homotopy connecting the maps $\ffe \colon N \to K$ to the set-valued map $\widehat \ff_0 \colon N \to 2^K$ with range $F = \widehat \ff_0(N)$, we proceed as follows.  As in the curve case, we construct its graph
$\barr\Gamma \subset \I \times  N \times K$
as the closure of the graph
$$\Gamma= \set{(\varepsilon ,\n,\ffe(\n)}{\n \in N,\ 0 < \varepsilon \leq \varepsilon _0}\subset (0,\varepsilon _0\,] \times  N \times K$$
 of the continuous map $F(\varepsilon ,\n) = \ffe(\n)$. Again, by the Closed Graph Theorem for set-valued functions, $\barr\Gamma$ is the graph of an upper hemicontinuous set-valued function $F \colon \I \times  N \to 2^K$, which,  for $0 < \varepsilon \leq \varepsilon _0$, defines the required upper hemicontinuous homotopy.
 
Finally, let us outline the proof in the general case.  To this end, we establish the following result by induction on the dimension $m$ of the face, using the preceding two-dimensional case as a model.

\begin{proposition} \label{F} Let $F \subset B$ be an $m$-dimensional face, and let $N \subset \Sdo$ be an $m$-dimensional submanifold satisfying the conditions of Proposition \ref{Nstar}. Then, the set-theoretic $\lim_{\varepsilon \to 0^+}\ffe(N) = F$.  Moreover, for $\varepsilon >0 $ sufficiently small, there is a continuous set valued homotopy from $\ffe \colon N \to K$ to the set-valued map $\widehat \ff_0\colon N \to 2^K$ whose range is the entire face: $\ff_0(N) = F$.
\end{proposition}

Referring back to the preceding argument for  two-dimensional polygonal faces, the key formulae \eq{Felx} work exactly as before, with $\subs \x k $ the vertices of $F$ and $\x_{k+1},\ldots, \x_{k+l}$ additional points in the configuration, if any, in $F\setminus \bsubs \x k$.  These in turn imply that, for $\varepsilon $ sufficiently small, the images $\ffe(N) $ lie in a $\O(\varepsilon )$ neighborhood of $F$, thus proving that $\lim_{\varepsilon \to 0^+} \ffe(N) \subset F$.

To prove that every point in $\x \in F$ is contained in the limit, we replace $N$ by the open submanifold $\widetilde N \subset N$ given in Proposition \ref{Nstar}. Again assume the contrary, that $\ffe(\widetilde N) \subset \Ue \setminus Z_{\x,\varepsilon}$, where $Z_{\x,\varepsilon} = Z_\x \capz \Ue$ with $Z_\x$ the affine subspace of dimension $d-m$ passing through $\x$ orthogonal to $F$.  According to the inductive hypothesis, its $(m-1)$-dimensional boundary component $L_i \subset \partial \widetilde N$ satisfies $\lim_{\varepsilon \to 0^+} \ffe(L_i) = G_i$, the corresponding $(m-1)$-dimensional subface of $F$,  through an upper hemicontinuous homotopy from $\ffe \colon L_i \to K$ to the set-valued map $\widehat \ff_0 \colon L_i \to 2^K$ with range $G_i = \widehat \ff_0(L_i) $.  We then, as in \eq{Gamma}, piece together these subface homotopies so as to construct an upper hemicontinuous homotopy from $\ffe \colon \partial \widetilde N \to K$ to the set-valued map $\widehat \ff_0 \colon \partial \widetilde N \to 2^K$ whose range is all of $\partial F = \widehat \ff_0(\partial \widetilde N)$.

The topological argument then proceeds in an identical manner, the only difference being that the open cover $\mathcal{V}$ is constructed as in \eq{Vo}, but now the intersections are homeomorphic to the contractible Cartesian product of a spherical sector of dimension $m$ with a ball of dimension $d-m$.  Further, we use the same commutative  diagram as in Figure  \ref{diagram} but with the first homology group $\mathcal{H}_1$ replaced by $\mathcal{H}_{m-1}$ throughout.  The resulting topological contradiction  proves that 
$$\lim_{\varepsilon \to 0^+} \ffe(\widetilde N) = \lim_{\varepsilon \to 0^+} \ffe(N) = F.$$ 
Finally, the construction of the corresponding upper hemicontinuous set-valued homotopy from $\ffe \colon N \to K$ to $\widehat \ff_0 \colon N \to 2^K$ with range $\widehat \ff_0(N) = F$ proceeds exactly as before.

The final step in the proof of the Main Theorem \ref{main_theorem} is to prove that $\lim_{\varepsilon \to 0^+} \ffe(\Sdo) = B$.  For this, we split up $B$ into its facets $B = F_1 \cupz \cdots \cupz F_k$.  For each $F_i$, by combining Lemma \ref{face_lemma} with the argument following \eq{Felx}, we deduce that $\lim_{\varepsilon \to 0^+} \ffe(S_{F_i}) \subset F_i$, and hence $\lim_{\varepsilon \to 0^+} \ffe(\Sdo) \subset B$. On the other hand, by the case $m = d-1$ of Proposition \ref F, there exists a $(d-1)$-dimensional submanifold $N_i \subset \Sdo$ such that $\lim_{\varepsilon \to 0^+} \ffe(N_i) = F_i$.  We conclude that $\lim_{\varepsilon \to 0^+} \ffe(\Sdo) = B$, as desired.  Moreover, we can similarly piece together the set-valued homotopies for each facet to find a set-valued homotopy from $\ffe \colon \Sdo \to K$ to the inverse Gauss map of the boundary, $\widehat \ff_0 \gamma _B^{-1}\colon \Sdo \to B$.
This, at last, completes our proof.

\appendix

\medskip

\section{Set-Valued Homology} \label{appB}

In this appendix, we review the basics of set-valued singular homology following Y.~Li, \cite{Li}. For simplicity, we will use $\mathbb{Q}$ as our ring of coefficients throughout.  

 Let $X$, $Y$ be connected normal Hausdorff topological spaces. Given a set-valued mapping $F\colon  X \to 2^Y$ and an open covering $\mathcal{U}$ of $Y$, we say that $F$ is \emph{compatible with $\mathcal{U}$} if and only if for any $x \in X$, there is 
some $U \in \mathcal{U}$ such that $F(x) \subset U$. Define
\Eq{CXYU}
$$\mathcal{C}(X,Y,\mathcal{U}) = \set{ \hskip-3pt F\colon  X \to 2^Y \hskip-3pt}{\hskip-10pt\text{ $F$ is an upper hemicontinuous mapping}\atop \text{\ with closed values compatible with $\>\mathcal{U}\>$} }.\sstrut{15}$$

Let 
$$\Delta_n = \set{x = \psubos xn \in \R^{n+1}}{x_i \geq 0, \ x_0 + x_1 + \cdotsx + x_n = 1}$$ denote the standard $n$-dimensional simplex. For $\rgo in$, let 
$$\varphi _i\psubos xn = (x_0, \ldots,x_{i-1},0,x_i,\ldots,x_n)$$ map the $(n-1)$-dimensional simplex $ \Delta_{n-1} $ to the $i$-th face $\displaystyle \Delta _n^{(i)} = \Delta _n \capz \{ x_i =  0\}$ of the $n$-dimensional simplex. 

Given an open cover $\mathcal{U}$ of $Y$, we define the \emph{$n$-th set-valued chain group} $\mathcal{C}_n(Y,\mathcal{U})$ to be the free abelian group\ generated by $\mathcal{C}(\Delta_n,Y,\mathcal{U})$ and call it the $n$-th set-valued 
chain group. We then define the boundary operator $\partial_n\colon  \mathcal{C}_n(Y,\mathcal{U}) \to \mathcal{C}_{n-1}(Y,\mathcal{U})$ by
\Eq{bdyop}
$$ \partial_n c_n = \sum_{i=0}^n \> (-1)^n\, c_n \comp \varphi _i. $$
Thus, $\partial_n \circ \partial_{n+1} = 0$, which is usually abbreviated by $\partial^2 = 0$. 

The $n$-th \emph{set-valued homology group} of $(Y,\mathcal{U})$ is then given by
\Eq{svHn}
$$ \mathcal{H}_n(Y, \mathcal{U}) = \Ker \partial_n / \Image \partial_{n+1}. $$
As noted by Li, \cite{Li}, if one is not careful when choosing the cover $\mathcal{U}$, all set-valued homology groups are trivial, and would thus be of no help establishing our desired topological result. To avoid this difficulty, Li imposes the \emph{contractible finite intersection property} on the cover $\mathcal{U}$.  This property requires that the intersection of any finite collection  of elements of the cover is either empty or contractible.

Since ordinary functions can be viewed as set-valued functions, there is a natural inclusion map $i$ from the $n$-th chain group $C_n(Y)$, as defined in the usual singular homology theory, to $\mathcal{C}_n(Y, \mathcal{U})$.  The inclusion is a chain map, and thus induces a group homomorphism
\Eq{isharp}
$$ i_{\sharp}\colon  H_n(Y) \longrightarrowz \mathcal{H}_n(Y, \mathcal{U}), $$
which, according to \cite[Theorem 11]{Li},  is actually an isomorphism.

Moreover, an upper hemicontinuous set-valued mapping $F\colon  X \to 2^Y$ with closed values  induces a chain map from $C_n(X)$ to $\mathcal{C}_n(Y, \mathcal{U})$, and thus induces a group homomorphism
\Eq{parenF}
$$ (F)_*\colon  H_n(X) \longrightarrowz \mathcal{H}_n(Y, \mathcal{U}). $$
In general, we will place parentheses around $(F)_*$ in order to distinguish the set-valued homology group homomorphism from the usual group homomorphism $f_* \colon H_n(X) \to H_n(Y)$ on the corresponding singular homology groups induced by a continuous (ordinary) function $f \colon X \to Y$.  Further results of Li, \cite{Li}, are quoted in the text as needed.

\bigskip

\acknowledge{We would like to thank Elias Saleeby, Dennis Sullivan, Daniele Tampieri, Paolo Emilio Ricci, and Kamal Khuri-Makdisi for their useful remarks and suggestions.}

\vskip8pt

\end{document}